\documentclass[11pt,a4paper,twoside]{article}

\usepackage{a4wide, amssymb, amsmath, amsthm, graphicx, xspace, enumerate}
\usepackage{algorithmic, algorithm}
\usepackage{ifpdf}

\usepackage{tikz}

\newtheorem{theorem}{Theorem}
\newtheorem{lemma}[theorem]{Lemma}
\newtheorem{corollary}[theorem]{Corollary}
\newtheorem{proposition}[theorem]{Proposition}
\newtheorem{conjecture}{Conjecture}

\newtheorem{invariant}{Invariant}

\newcommand{\inst}[1]{$^{#1}$}

\newcommand{\cone}[1]{\ensuremath{C_1(#1)}}
\newcommand{\czero}[1]{\ensuremath{C_0(#1)}}
\newcommand{\ctwo}[1]{\ensuremath{C_2(#1)}}
\newcommand{\ceil}[1]{\lceil #1 \rceil}
\newcommand{\floor}[1]{\lfloor #1 \rfloor}

\newcommand{\ignore}[1]{}

\ifpdf
\newcommand{\figext}{pdf}
\else
\newcommand{\figext}{eps}
\fi

\begin{document}
\date{}
\title{A Planar Linear Arboricity Conjecture\thanks{Supported in part by bilateral project BI-PL/08-09-008. M. Cygan and \L. Kowalik were supported in part by Polish Ministry of Science and 
		Higher Education grant N206 355636}
\author{Marek Cygan\inst{1},
	{\L}ukasz Kowalik\inst{1},
  	Borut Lu\v{z}ar\inst{2}\thanks{Operation part financed by the European Union, European Social Fund.}}}

\maketitle
\begin{center}
{\footnotesize
\inst{1} Institute of Informatics, Warsaw University, \\
         Banacha 2, 02-097 Warsaw, Poland\\
         \texttt{\{cygan,kowalik\}@mimuw.edu.pl}\\
\inst{2} Institute of Mathematics, Physics, and Mechanics,\\
         Jadranska~19, 1111 Ljubljana, Slovenia\\
         \texttt{borut.luzar@gmail.com}
} \end{center}

\begin{abstract}
The linear arboricity ${\rm la}(G)$ of a graph $G$ is the minimum number of linear forests that partition
the edges of $G$. In 1984, Akiyama et al.~\cite{AEH} stated the Linear Arboricity Conjecture (LAC), that the linear arboricity
of any simple graph of maximum degree $\Delta$ is either $\big \lceil \tfrac{\Delta}{2} \big \rceil$ or 
$\big \lceil \tfrac{\Delta+1}{2} \big \rceil$. In~\cite{Wu,WW} it was proven that LAC holds for all planar graphs. 

LAC implies that for $\Delta$ odd, ${\rm la}(G)=\big \lceil \tfrac{\Delta}{2} \big \rceil$.
We conjecture that for planar graphs this equality is true also for any even $\Delta \ge 6$.
In this paper we show that it is true for any even $\Delta \ge 10$, leaving open
only the cases $\Delta=6, 8$.

We present also an $O(n\log n)$ algorithm for partitioning a planar graph into $\max\{{\rm la}(G), 5\}$ linear forests,
which is optimal when $\Delta \ge 9$.
\end{abstract}


%
%
%
%

\section{Introduction}

In this paper we consider only undirected and simple graphs. A {\em linear forest} is a
forest in which every connected component is a path. The linear arboricity
${\rm la}(G)$ of a graph $G$ is the minimum number of linear forests in $G$, whose union
is the set of all edges of $G$. This one of the most natural graph covering notions was introduced by Harary~\cite{harary} in 1970. 
Note that for any graph of maximum degree $\Delta$ one needs at least $\ceil{\tfrac{\Delta}{2}}$ linear forests to cover all the edges.
If $\Delta$ is even and the graph is regular, $\ceil{\tfrac{\Delta}{2}}$ forests do not suffice, for otherwise every vertex in every forest
has degree $2$, a contradiction. Hence, for any $\Delta$-regular graph $G$, we have ${\rm la}(G)\ge\ceil{\tfrac{\Delta+1}{2}}$. 
Akiyama, Exoo and Harary conjectured that this bound is always tight, i.e.\ for any regular graph, ${\rm la}(G)=\ceil{\tfrac{\Delta+1}{2}}$.
It is easy to see (check e.g.~\cite{alon_ijm}) that this conjecture is equivalent to

\begin{conjecture}
\label{lac}
For any graph $G$,  $\ceil{\tfrac{\Delta}{2}}\le{\rm la}(G)\le\ceil{\tfrac{\Delta+1}{2}}$.
\end{conjecture}

We note that Conjecture~\ref{lac} resembles Vizing Theorem, and indeed, for odd $\Delta$ it can be treated as a generalization of Vizing Theorem (just color each linear forest with two new colors and get a $(\Delta+1)$-edge-coloring). However, despite many efforts, the conjecture is still open and the
best known upper general bound is ${\rm la}(G)=\tfrac{d}{2} + O(d^{2/3}(\log d)^{1/3})$, due to Alon, Teague and Wormald~\cite{ATW}. 
Conjecture~\ref{lac} was proved only in several special cases: for $\Delta=3,4$ in~\cite{AEH,akiyama2}, for $\Delta=5,6,8$ in~\cite{EP}, $\Delta=10$~\cite{Guldan}, for bipartite graphs~\cite{AEH} to mention a few. Finally, it was shown for planar graphs (for $\Delta\ne 7$ by Jian-Liang Wu~\cite{Wu} and for $\Delta=7$ by Jian-Liang Wu and Yu-Wen Wu). 

In this paper we focus on planar graphs. Note that for all odd $\Delta$, we have $\ceil{\tfrac{\Delta}{2}}=\ceil{\tfrac{\Delta+1}{2}}$.
Hence, the linear arboricity is $\ceil{\tfrac{\Delta}{2}}$ for planar graphs of odd $\Delta$. 
Moreover, in~\cite{Wu}, Wu showed that this is also true for $\Delta \ge 13$.
In this paper we state the following conjecture.

\begin{conjecture}
\label{plac}
For any planar graph $G$ of maximum degree $\Delta\ge 6$, we have ${\rm la}(G)=\ceil{\tfrac{\Delta}{2}}$.
\end{conjecture}

It is easy to see that the above equality does not hold for $\Delta=2,4$, since there are $2$- and $4$-regular planar graphs: e.g. the $3$-cycle and
the octahedron. Interestingly, if the equality holds for $\Delta=6$, then any planar graph with maximum degree $6$ is $6$-edge-colorable (just edge-color each of the three linear forests in two new colors). 
Hence, Conjecture~\ref{plac} implies the Vizing Planar Graph Conjecture~\cite{Vizing65} (as currently it is open only for $\Delta=6$):

\begin{conjecture}[Vizing Planar Graph Conjecture]
\label{vizing}
For any planar graph $G$ of maximum degree $\Delta\ge 6$, we have $\chi'(G)=\Delta$.
\end{conjecture}

The main result of this paper is the following theorem.

\begin{theorem}
	\label{thm_main}
	For every planar graph $G$ of maximum degree $\Delta(G) \ge 9$, we have ${\rm la}(G)=\ceil{\tfrac{\Delta}{2}}$.
\end{theorem}

We note that Wu, Hou and Sun~\cite{WHS} verified Conjecture~\ref{plac} for planar graphs without $4$-cycles.
For $\Delta\ge 7$ it is also known to be true for planar graphs without $3$-cycles~\cite{Wu} and without $5$-cycles~\cite{WHL}.

\paragraph{Computational Complexity Perspective.} 
Consider the following decision problem. Given graph $G$ and a number $k$, determine whether ${\rm la}(G)=k$.
Peroche~\cite{peroche} showed that this problem is NP-complete even for graphs of maximum degree $\Delta=4$.
Our result settles the complexity of the decision problem for planar graphs of maximum degree $\Delta(G) \ge 9$.
The discussion above implies that for planar graphs the decision problem is trivial also when $\Delta$ is odd.
When $\Delta=2$ the problem is in $P$, as then the algorithm just checks whether there is a cycle in $G$.
Hence, the remaining cases are $\Delta=4,6,8$. Conjecture~\ref{plac}, if true, excludes also cases $\Delta=6,8$.
We conjecture that the only remaining case $\Delta=4$ is NP-complete.

\begin{conjecture}
\label{con:np}
It is NP-complete to determine whether a given planar graph of maximum degree $4$ has linear arboricity $2$.
\end{conjecture}

Finally, even when one knows the linear arboricity of a given graph (or a bound on it) the question arises how fast one can find
the corresponding collection of linear forests. We show the following result.

\begin{theorem}
	\label{thm_algo}
	For every $n$-vertex planar graph $G$ of maximum degree $\Delta$, one can find a cover of its edges by
	$\max\{{\rm la}(G),5\}$ linear forests in $O(n\log n)$ time.
\end{theorem}

\paragraph{Preliminaries.} 
        For any graph $G$, $V(G)$ and $E(G)$ denote the sets of vertices and edges of $G$, respectively. If $G$ is plane, we also denote the set of its
	faces by $F(G)$.

	We call a vertex of degree $k$, at least $k$, and at most $k$, a $k$-vertex, $(\ge  k)$-vertex, and $(\le  k)$-vertex, respectively.
	Moreover, a neighbor $u$ of a vertex $v$ is called a $k$-neighbor, where $k = \deg(u)$.
	In a plane graph, {\em length} of a face $f$, denoted as $\ell(f)$ is the number of edges incident to $f$.
	Analogously to vertices, a face of length $k$, at least $k$, and at most $k$, is called a $k$-face, $(\ge  k)$-face, and $(\le  k)$-face, respectively.
	
	Throughout the paper it will be convenient for us to treat partitions into linear forests as a kind of edge-colorings.
	A {\em $k$-linear coloring} of a graph $G=(V,E)$ is a function $C:E\rightarrow \{1, \ldots, k\}$ such that for $i=1,\ldots,k$, 
	the set of edges $C^{-1}(i)$ is a linear forest. We call an edge colored by $a$, an $a$-edge.	
	
	Following the notation in~\cite{Wu}, by $\czero{v}$, we denote the set of colors that are not used on edges incident to vertex $v$, and
	$\cone{v}$ is the set of colors which are assigned to exactly one edge incident to $v$. Finally, $\ctwo{v}=\{1,\ldots,k\}\setminus(\czero{v}\cup\cone{v})$.
	We call the color which is in $\czero{v} \cup \cone{v}$ a \textit{free} color at $v$. 
	A  path, where all edges have the same color $a$ is called a {\em monochromatic path} or an {\em $a$-path}.
	A {\em $(uv,a)$-path} is an $a$-path with endpoints $u$ and $v$.
	
	We use also the Iverson's notation, i.e.\ $[\alpha]$ is 1 when $\alpha$ is true and 0 otherwise.
	
%
%
%
%

\section{Proof of Theorem~\ref{thm_main}}

It will be convenient for us to prove the following generalization of Theorem~\ref{thm_main}:

\begin{proposition}
\label{prop:main}
Any simple planar graph of maximum degree $\Delta$ has a linear coloring in $\max\{\ceil{\tfrac{\Delta}{2}}, 5\}$ colors.
\end{proposition}

Our plan for proving Proposition~\ref{prop:main} is as follows. 
First we are going to show some structural properties of a graph which is {\em not} $k$-linear colorable and,
subject to this condition, is minimal, i.e.\ has as few edges as possible.
In this phase we do not use planarity, and our results apply to general graphs.
Next, using the so-called discharging method (a powerful technique developed for proving the four color theorem)
we will show that when $k=\max\{\ceil{\tfrac{\Delta}{2}}, 5\}$, there is no planar graph with the obtained structural properties. 

\subsection{Structure of a minimal graph $G$ with $\Delta(G)\le 2k$ and ${\rm la}(G) > k$}
\label{sec:structure}

In this section we fix a number $k$ and we assume that $G$ is a graph of maximum degree at most $2k$ which is {\em not} $k$-linear colorable
and, among such graphs, $G$ has minimal number of edges.
The following Lemma appears in many previous works, e.g.\ in~\cite{WHL}, but we give the proof for completeness.

\begin{lemma}
	\label{lem_edge}
	For every edge $uv$ of $G$, $d(u) + d(v) \ge 2k + 2$.	
\end{lemma}

\begin{proof}
	Suppose for a contradiction that $uv$ is an edge in $G$, and $d(u) + d(v) < 2k + 2$.
	By minimality of $G$, there exists a $k$-linear coloring of $G'=G - uv$. Note
	that the degree of vertices $u$ and $v$ in $G'$ is one less than in $G$, hence
	$d_{G'}(u) + d_{G'}(v) < 2k$. 
	So there exists at least one color $c$ which is either an element of $\czero{u}$ and free at $v$, or 
	an element of $\czero{v}$ and free at $u$. 
	It follows that we can extend the $k$-linear coloring of $G'$ to $G$, a contradiction.
\end{proof}

\begin{lemma}
\label{lem:2vertex}
  $G$ does not contain a $2$-vertex $v$ such that the two neighbors $u$ and $z$ are not adjacent.
\end{lemma}

\begin{proof}
Since there is no edge $uz$ we can create a simple graph $G'=(G - v) \cup uz$.
Because of the minimality of $G$ there exists a $k$-linear coloring $C$ of $G'$.
Let $a = C(uz)$. We can put $C(uv)=C(vz)=a$ obtaining a $k$-linear coloring of $G$,
a contradiction.
\end{proof}

\begin{lemma}
  \label{lem:atmosttwo2vertices}
	Every vertex has at most one adjacent $2$-vertex.
\end{lemma}

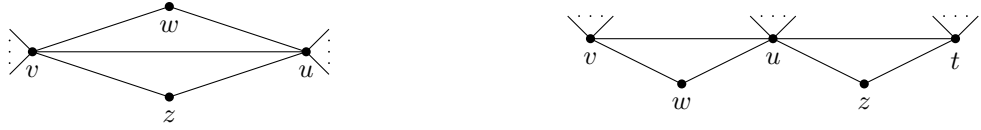
\begin{figure}[ht]
\centering
\begin{tabular}{cc}

\begin{minipage}{2.5in}
\centering

{\small{
\begin{tikzpicture}[scale=0.6]
  \begin{scope}[shift={(10.0, 0.0)}]

  \foreach \x/\y/\z in {0/0/$v$, 3/1/$w$, 6/0/$u$, 3/-1/$z$}
  {
    \fill (\x,\y) circle (0.1);
    \draw[below] (\x,\y-0.1) node {\z};
  }

  \draw (0,0) -- (3,1) -- (6,0) -- (3,-1) -- (0,0);
  \draw (0,0) -- (6,0);

  \foreach \y in {-0.25, 0, 0.25}
  {
    \fill (-0.5, \y) circle (0.02);
  }
  \draw (-0.5,0.5) -- (0,0) -- (-0.5,-0.5);

  \foreach \y in {-0.25, 0, 0.25}
  {
    \fill (6.5, \y) circle (0.02);
  }
  \draw (6.5,0.5) -- (6,0) -- (6.5,-0.5);

  \end{scope}
\end{tikzpicture}
}}
\end{minipage}

&

\begin{minipage}{2.5in}
\centering
{\small{
\begin{tikzpicture}[scale=0.6]

  \foreach \x in {-4, 0, 4}
  {
    \fill (\x-0.25, 0.5) circle (0.02);
    \fill (\x,0.5) circle (0.02);
    \fill (\x+0.25, 0.5) circle (0.02);
    \draw (\x-0.5,0.5) -- (\x,0) -- (\x+0.5,0.5);
  }

  \foreach \x/\y/\z in {-4/0/$v$, 0/0/$u$, 4/0/$t$, -2/-1/$w$, 2/-1/$z$}
  {
    \fill (\x,\y) circle (0.1);
    \draw[below] (\x,\y-0.1) node {\z};
  }

  \draw (-4,0) -- (4,0) -- (2,-1) -- (0,0) -- (-2,-1) -- (-4,0);
\end{tikzpicture}
}}
\end{minipage}

\\

(A) $2$-vertices have two common neighbors & 

(B) $2$-vertices have only one common neighbor
\\

\end{tabular}

\caption{Cases $A$ and $B$ in the proof of Lemma~\ref{lem:atmosttwo2vertices}}
\label{fig:two2vertices}
\end{figure}

\begin{proof}
  Assume for a contradiction that a vertex $u$ has at least two neighbors of degree two namely $w$ and $z$.
  By Lemma~\ref{lem:2vertex}, neighbors of every $2$-vertex are adjacent.
  Consider the configuration and the labeling of vertices as in Fig.~\ref{fig:two2vertices}.

  $Case\ A.$ \textit{The vertices $w$ and $z$ have two common neighbors}.
  Since $G$ is minimal, there exists a $k$-linear coloring $C$ of $G'=G-uz$.
  Let $a$ be a free color at $u$. We may assume that $C(vz)=a$ and there exists
  a $(uv,a)$-path, otherwise we can put $C(uz)=a$ and we are done.
  Consider edges $vw$ and $wu$. If at least one of them has color $a$ it means that
  both have color $a$ since there is a $(uv,a)$-path. 
  In such a case ($C(vw)=C(wu)=a$) we can take color $b=C(uv)\not =a$ and recolor
  $C(vw)=C(wu)=b$ and $C(uv)=a$. 
  Thus we can assume that $C(vw)=b\not =a$ and $C(wu)=c\not = a$ therefore it
  suffices to recolor $C(uw)=a$ and put $C(uz)=c$.

\begin{figure}[ht]
\begin{center}
{\small{
\begin{tikzpicture}[scale=0.55]

  \foreach \x/\y/\z in {0/0/$v$, 3/1/$w$, 6/0/$u$, 3/-1/$z$}
  {
    \fill (\x,\y) circle (0.1);
    \draw[below] (\x,\y-0.1) node {\z};
  }

  \draw (0,0) -- (3,1) -- (6,0) -- (0,0);
  \draw (0,0) -- (3,-1);
  \draw[dashed](3,-1) -- (6,0);

  \foreach \y in {-0.25, 0, 0.25}
  {
    \fill (-0.5, \y) circle (0.02);
  }
  \draw (-0.5,0.5) -- (0,0) -- (-0.5,-0.5);

  \foreach \y in {-0.25, 0, 0.25}
  {
    \fill (6.5, \y) circle (0.02);
  }
  \draw (6.5,0.5) -- (6,0) -- (6.5,-0.5);

  \draw[below left] (1.5,-0.5) node {$a$};
  \draw[below right] (4.5,-0.5) node {$a$};

  \draw[above left] (1.5,0.5) node {$b$};
  \draw[above right] (4.5,0.5) node {$c$};

  \draw[dashed,thick,<-] (10,1) arc (65:115:3*1.414);

  \begin{scope}[shift={(10.0, 0.0)}]
    \foreach \x/\y/\z in {0/0/$v$, 3/1/$w$, 6/0/$u$, 3/-1/$z$}
    {
      \fill (\x,\y) circle (0.1);
      \draw[below] (\x,\y-0.1) node {\z};
    }

    \draw (0,0) -- (3,1) -- (6,0) -- (0,0);
    \draw (0,0) -- (3,-1);
    \draw (3,-1) -- (6,0);

    \foreach \y in {-0.25, 0, 0.25}
    {
      \fill (-0.5, \y) circle (0.02);
    }
    \draw (-0.5,0.5) -- (0,0) -- (-0.5,-0.5);

    \foreach \y in {-0.25, 0, 0.25}
    {
      \fill (6.5, \y) circle (0.02);
    }
    \draw (6.5,0.5) -- (6,0) -- (6.5,-0.5);

    \draw[below left] (1.5,-0.5) node {$a$};
    \draw[below right] (4.5,-0.5) node {$c$};

    \draw[above left] (1.5,0.5) node {$b$};
    \draw[above right] (4.5,0.5) node {$a$};
  \end{scope}

\end{tikzpicture}
}}
\end{center}
\end{figure}

  $Case\ B.$ \textit{Vertices $w$ and $z$ have precisely one common neighbor}.
  Since $G$ is minimal, there exists a $k$-linear coloring $C$ of $G'=G-uz$.
  Let $a$ be the only free color at $u$. We may assume that $C(tz)=a$ and there exists
  a $(tu,a)$-path since otherwise we can put $C(uz)=a$ and we are done.
  Now let us take into consideration the edge $uw$.

  $Case\ B.1.$ $C(uw)=a$. Since in $G'$ the only $a$-edge incident to the vertex $u$ is the edge $uw$
  we have $C(wv)=a$ because there exists a $(tu,a)$-path in $G'$. 
  Let $C(uv)=b$. 
  We know that $b\not =a$ since $a$ is free at $u$.
  In this case we recolor $C(vw)=b$, $C(uv)=a$ and $C(uz)=b$.
  It is easy to see that we do not introduce any monochromatic cycle.

\begin{figure}[hptb]
\begin{center}
{\small{
\begin{tikzpicture}[scale=0.55]

  \foreach \x in {-4, 0, 4}
  {
    \fill (\x-0.25, 0.5) circle (0.02);
    \fill (\x,0.5) circle (0.02);
    \fill (\x+0.25, 0.5) circle (0.02);
    \draw (\x-0.5,0.5) -- (\x,0) -- (\x+0.5,0.5);
  }

  \foreach \x/\y/\z in {-4/0/$v$, 0/0/$u$, 4/0/$t$, -2/-1/$w$, 2/-1/$z$}
  {
    \fill (\x,\y) circle (0.1);
    \draw[below] (\x,\y-0.1) node {\z};
  }

  \draw (-4,0) -- (4,0) -- (2,-1);
  \draw (0,0) -- (-2,-1) -- (-4,0);
  \draw[dashed] (0,0) -- (2,-1);

  \draw[below left] (-3,-0.5) node {$a$};
  \draw[below right] (-1,-0.5) node {$a$};
  \draw[above] (-2,0) node {$b$};
  \draw[below left] (1,-0.5) node {$a$};
  \draw[below right] (3,-0.5) node {$a$};

  \draw[dashed,thick,<-] (7,1) arc (65:115:3*1.414);

  \begin{scope}[shift={(10.0, 0.0)}]
    \foreach \x in {-4, 0, 4}
    {
      \fill (\x-0.25, 0.5) circle (0.02);
      \fill (\x,0.5) circle (0.02);
      \fill (\x+0.25, 0.5) circle (0.02);
      \draw (\x-0.5,0.5) -- (\x,0) -- (\x+0.5,0.5);
    }

    \foreach \x/\y/\z in {-4/0/$v$, 0/0/$u$, 4/0/$t$, -2/-1/$w$, 2/-1/$z$}
    {
      \fill (\x,\y) circle (0.1);
      \draw[below] (\x,\y-0.1) node {\z};
    }

    \draw (-4,0) -- (4,0) -- (2,-1);
    \draw (0,0) -- (-2,-1) -- (-4,0);
    \draw (0,0) -- (2,-1);

    \draw[below left] (-3,-0.5) node {$b$};
    \draw[below right] (-1,-0.5) node {$a$};
    \draw[above] (-2,0) node {$a$};
    \draw[below left] (1,-0.5) node {$b$};
    \draw[below right] (3,-0.5) node {$a$};
  \end{scope}

\end{tikzpicture}
}}
\end{center}
\end{figure}

  $Case\ B.2.$ $C(uw)=b\not =a$. In this case we only recolor $C(uw)=a$ and put $C(uz)=b$.
  Even if $C(vw)=a$ we do not introduce a $a$-cycle since then a $(uv,a)$-path can not exist.
\begin{center}
{\small{
\begin{tikzpicture}[scale=0.55]

  \foreach \x in {-4, 0, 4}
  {
    \fill (\x-0.25, 0.5) circle (0.02);
    \fill (\x,0.5) circle (0.02);
    \fill (\x+0.25, 0.5) circle (0.02);
    \draw (\x-0.5,0.5) -- (\x,0) -- (\x+0.5,0.5);
  }

  \foreach \x/\y/\z in {-4/0/$v$, 0/0/$u$, 4/0/$t$, -2/-1/$w$, 2/-1/$z$}
  {
    \fill (\x,\y) circle (0.1);
    \draw[below] (\x,\y-0.1) node {\z};
  }

  \draw (-4,0) -- (4,0) -- (2,-1);
  \draw (0,0) -- (-2,-1) -- (-4,0);
  \draw[dashed] (0,0) -- (2,-1);

  \draw[below right] (-1,-0.5) node {$b$};
  \draw[below left] (1,-0.5) node {$a$};
  \draw[below right] (3,-0.5) node {$a$};
  \draw[dashed,thick,<-] (7,1) arc (65:115:3*1.414);

  \begin{scope}[shift={(10.0, 0.0)}]
    \foreach \x in {-4, 0, 4}
    {
      \fill (\x-0.25, 0.5) circle (0.02);
      \fill (\x,0.5) circle (0.02);
      \fill (\x+0.25, 0.5) circle (0.02);
      \draw (\x-0.5,0.5) -- (\x,0) -- (\x+0.5,0.5);
    }

    \foreach \x/\y/\z in {-4/0/$v$, 0/0/$u$, 4/0/$t$, -2/-1/$w$, 2/-1/$z$}
    {
      \fill (\x,\y) circle (0.1);
      \draw[below] (\x,\y-0.1) node {\z};
    }

    \draw (-4,0) -- (4,0) -- (2,-1);
    \draw (0,0) -- (-2,-1) -- (-4,0);
    \draw (0,0) -- (2,-1);

    \draw[below right] (-1,-0.5) node {$a$};
    \draw[below left] (1,-0.5) node {$b$};
    \draw[below right] (3,-0.5) node {$a$};
  \end{scope}
\end{tikzpicture}
}}
\end{center}
\end{proof}

\begin{lemma}
\label{lem:two-pairs}
	$G$ does not contain a $3$-vertex with precisely two pairs of adjacent neighbors.
\end{lemma}

\begin{figure}[ht]
	\centering
        \includegraphics{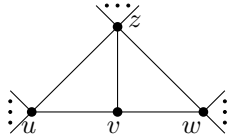}
	\caption{A $3$-vertex with precisely two pairs of adjacent neighbors.}
	\label{fig_cubtri}
\end{figure}

\begin{proof}
	Assume the configuration described in the claim exists and consider the labeling of vertices as in Fig.~\ref{fig_cubtri}, in particular $uw\not\in E(G)$.
	Let $G'$ be the graph obtained by removing $v$ and adding the edge
	$uv$ to $G$. By the minimality of $G$, there exists a $k$-linear coloring $C$ of $G'$. Let $a$, $b$, and $c$ be the colors
	of the edges $uz$, $zw$, and $uw$, respectively. We show that $C$ can be extended to $G$ as follows.
	
	First, color the edges in $E(G)\cap E(G')$ as in the coloring $C$. Then only the edges $uv$, $vw$, and $vz$ remain non-colored. 
	Let $d$ be a color free at $z$. 
	If $d \neq c$ we color $vz$ by $d$ and both $vu$, $vw$ by $c$. 
	It is easy to see that this is a proper $k$-linear coloring of $G$ (note that there is no $(uw,c)$-path in $G$, since then there is a $c$-cycle 
	in the coloring $C$).
	So, we may assume that $c$ is free at $z$.
	Now observe that if in the partial coloring of $G$ there are both a $(uz,c)$-path and $(wz,c)$-path then $G'$ contains a $c$-cycle.
	By symmetry, we can assume that there is no $(uz,c)$-path. Then, in particular, $c\ne a$.
		We color the edges $uv$ and $vz$ with $a$, $vw$, $uz$ with $c$, and $wz$ remains colored with $b$.
	Thus we obtained a $k$-linear coloring of $G$, a contradiction.
\end{proof}

\begin{figure}[h]
\begin{center}
{\small{
\begin{tikzpicture}[scale=0.5]

  \foreach \x/\y/\z in {0/0/$w$, 2/2/$u$, 0/4/$z$}
  {
    \fill (\x,\y) circle (0.1);
    \draw[right] (\x,\y) node {\z};
  }

  \fill (-2,2) circle (0.1);
  \draw[above] (-2,2) node {$v$};

  \draw (0,0) -- (2,2) -- (0,4) -- (-2,2) -- (0,0);
  \draw (0,0) -- (0,4);

  \draw (-2,2) -- (-2.5,2);
  \draw (0,0) -- (-0.5,-0.5);
  \draw (0,0) -- (0.5,-0.5);
  \fill (-0.25, -0.5) circle (0.02);
  \fill (0,-0.5) circle (0.02);
  \fill (0.25, -0.5) circle (0.02);

  \draw (0,4) -- (-0.5,4.5);
  \draw (0,4) -- (0.5,4.5);
  \fill (-0.25, 4.5) circle (0.02);
  \fill (0,4.5) circle (0.02);
  \fill (0.25, 4.5) circle (0.02);

\end{tikzpicture}
}}
\caption{Configuration from Lemma~\ref{lem:23triangles}}
\label{fig:23triangles}
\end{center}
\end{figure}
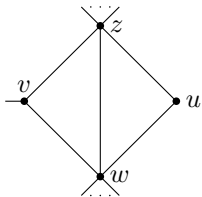

\begin{lemma}
\label{lem:23triangles}
$G$ does not contain the configuration of Fig.~\ref{fig:23triangles}, i.e.\ a 4-cycle $vzuw$ with a chord $zw$ such that $\deg(v)=3$ and $\deg(u)=2$.
\end{lemma}

\begin{proof}
Consider the configuration and the labeling of vertices as in Fig.~\ref{fig:23triangles}.
Since $G$ is a minimal counterexample, there exists a $k$-linear coloring $C$ of $G'=G-uz$.
Let $a$ be a free color at $z$. We may assume that $C(uw)=a$ and there exists
a $(wz,a)$-path since otherwise we can put $C(uz)=a$ and we are done.
Now let us take into consideration the edge $wz$.

  $Case\ 1.$ $C(wz)=a$. then since $a$ is free at the vertex $v$ we know that $C(vz)=b\not =a$
  and $C(vw)=c\not =a$ since the vertex $w$ already has two incident edges of color $a$.
  In this case we recolor $C(vz)=a$ and put $C(uz)=b$ and we do not introduce any monochromatic cycle.

\begin{figure}[ht]
\begin{center}
{\small{
\begin{tikzpicture}[scale=0.5]

  \foreach \x/\y/\z in {0/0/$w$, 2/2/$u$, 0/4/$z$}
  {
    \fill (\x,\y) circle (0.1);
    \draw[right] (\x,\y) node {\z};
  }

  \fill (-2,2) circle (0.1);
  \draw[above] (-2,2) node {$v$};

  \draw (0,0) -- (2,2);
  \draw (0,4) -- (-2,2) -- (0,0);
  \draw[dashed] (2,2) -- (0,4);
  \draw (0,0) -- (0,4);

  \draw (-2,2) -- (-2.5,2);
  \draw (0,0) -- (-0.5,-0.5);
  \draw (0,0) -- (0.5,-0.5);
  \fill (-0.25, -0.5) circle (0.02);
  \fill (0,-0.5) circle (0.02);
  \fill (0.25, -0.5) circle (0.02);

  \draw (0,4) -- (-0.5,4.5);
  \draw (0,4) -- (0.5,4.5);
  \fill (-0.25, 4.5) circle (0.02);
  \fill (0,4.5) circle (0.02);
  \fill (0.25, 4.5) circle (0.02);

  \draw[below right] (1,1) node {$a$};
  \draw[above right] (1,3) node {$a$};
  \draw[right] (0,2) node {$a$};
  \draw[above left] (-1,3) node {$b$};
  \draw[below left] (-1,1) node {$c$};

  \draw[very thick, dashed, ->] (3.5, 2.0) -- (4.5, 2.0);

  \begin{scope}[shift={(8.0, 0.0)}]
    \foreach \x/\y/\z in {0/0/$w$, 2/2/$u$, 0/4/$z$}
    {
      \fill (\x,\y) circle (0.1);
      \draw[right] (\x,\y) node {\z};
    }

    \fill (-2,2) circle (0.1);
    \draw[above] (-2,2) node {$v$};

    \draw (0,0) -- (2,2);
    \draw (0,4) -- (-2,2) -- (0,0);
    \draw (0,0) -- (0,4);

    \draw (-2,2) -- (-2.5,2);
    \draw (0,0) -- (-0.5,-0.5);
    \draw (0,0) -- (0.5,-0.5);
    \fill (-0.25, -0.5) circle (0.02);
    \fill (0,-0.5) circle (0.02);
    \fill (0.25, -0.5) circle (0.02);

    \draw (0,4) -- (-0.5,4.5);
    \draw (0,4) -- (0.5,4.5);
    \fill (-0.25, 4.5) circle (0.02);
    \fill (0,4.5) circle (0.02);
    \fill (0.25, 4.5) circle (0.02);
    \draw (2,2) -- (0,4);

    \draw[above right] (1,3) node {$b$};
    \draw[below right] (1,1) node {$a$};
    \draw[right] (0,2) node {$a$};
    \draw[above left] (-1,3) node {$a$};
    \draw[below left] (-1,1) node {$c$};
  \end{scope} 
\end{tikzpicture}
}}
\caption{Case $1$ of the proof of Lemma~\ref{lem:23triangles}.}
\end{center}
\end{figure}
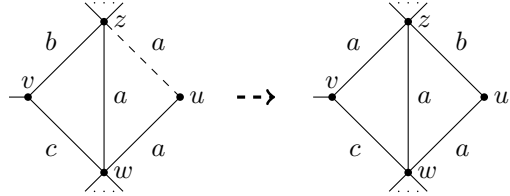

  $Case\ 2.$ $C(wz)=b\not =a$. Now we consider the edges $vz$ and $vw$:

  $Case\ 2.1$ $C(wz)=b\not =a$, $C(vz)=c\not =a$, $C(vw)=d\not =a$. In this case
  we recolor $C(vz)=a$ and put $C(uz)=c$. We do not introduce any $a$-cycle 
  because there is a ($wz$,$a$)-path which means that there is no ($vw$,$a$)-path
  since $w$ has only one adjacent $a$-edge outside of the configuration.

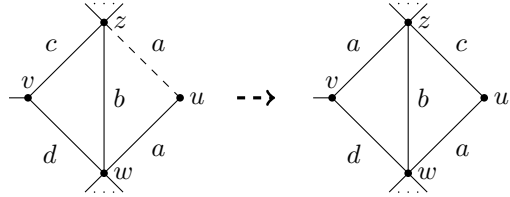
\begin{figure}[h]
\begin{center}
{\small{
\begin{tikzpicture}[scale=0.5]

  \foreach \x/\y/\z in {0/0/$w$, 2/2/$u$, 0/4/$z$}
  {
    \fill (\x,\y) circle (0.1);
    \draw[right] (\x,\y) node {\z};
  }

  \fill (-2,2) circle (0.1);
  \draw[above] (-2,2) node {$v$};

  \draw (0,0) -- (2,2);
  \draw (0,4) -- (-2,2) -- (0,0);
  \draw[dashed] (2,2) -- (0,4);
  \draw (0,0) -- (0,4);

  \draw (-2,2) -- (-2.5,2);
  \draw (0,0) -- (-0.5,-0.5);
  \draw (0,0) -- (0.5,-0.5);
  \fill (-0.25, -0.5) circle (0.02);
  \fill (0,-0.5) circle (0.02);
  \fill (0.25, -0.5) circle (0.02);

  \draw (0,4) -- (-0.5,4.5);
  \draw (0,4) -- (0.5,4.5);
  \fill (-0.25, 4.5) circle (0.02);
  \fill (0,4.5) circle (0.02);
  \fill (0.25, 4.5) circle (0.02);

  \draw[below right] (1,1) node {$a$};
  \draw[above right] (1,3) node {$a$};
  \draw[right] (0,2) node {$b$};
  \draw[above left] (-1,3) node {$c$};
  \draw[below left] (-1,1) node {$d$};

  \draw[very thick, dashed, ->] (3.5, 2.0) -- (4.5, 2.0);

  \begin{scope}[shift={(8.0, 0.0)}]
    \foreach \x/\y/\z in {0/0/$w$, 2/2/$u$, 0/4/$z$}
    {
      \fill (\x,\y) circle (0.1);
      \draw[right] (\x,\y) node {\z};
    }

    \fill (-2,2) circle (0.1);
    \draw[above] (-2,2) node {$v$};

    \draw (0,0) -- (2,2);
    \draw (0,4) -- (-2,2) -- (0,0);
    \draw (0,0) -- (0,4);

    \draw (-2,2) -- (-2.5,2);
    \draw (0,0) -- (-0.5,-0.5);
    \draw (0,0) -- (0.5,-0.5);
    \fill (-0.25, -0.5) circle (0.02);
    \fill (0,-0.5) circle (0.02);
    \fill (0.25, -0.5) circle (0.02);

    \draw (0,4) -- (-0.5,4.5);
    \draw (0,4) -- (0.5,4.5);
    \fill (-0.25, 4.5) circle (0.02);
    \fill (0,4.5) circle (0.02);
    \fill (0.25, 4.5) circle (0.02);
    \draw (2,2) -- (0,4);

    \draw[above right] (1,3) node {$c$};
    \draw[below right] (1,1) node {$a$};
    \draw[right] (0,2) node {$b$};
    \draw[above left] (-1,3) node {$a$};
    \draw[below left] (-1,1) node {$d$};
  \end{scope} 
\end{tikzpicture}
}}
\caption{Case $2.1$ of the proof of Lemma~\ref{lem:23triangles}.}
\end{center}
\end{figure}

  $Case\ 2.2$ $C(wz)=b\not =a$, $C(vz)=c\not =a$, $C(vw)=a$. In this case
  we recolor $C(wz)=a$, $C(vw)=b$ and put $C(uz)=b$.
  Since there was a ($wz$,$a$)-path it means that the only outside edge of the vertex
  $v$ has color $a$ thus even if $b=c$ we do not introduce any monochromatic cycle.

\begin{figure}[h]
\begin{center}
{\small{
\begin{tikzpicture}[scale=0.5]

  \foreach \x/\y/\z in {0/0/$w$, 2/2/$u$, 0/4/$z$}
  {
    \fill (\x,\y) circle (0.1);
    \draw[right] (\x,\y) node {\z};
  }

  \fill (-2,2) circle (0.1);
  \draw[above] (-2,2) node {$v$};

  \draw (0,0) -- (2,2);
  \draw (0,4) -- (-2,2) -- (0,0);
  \draw[dashed] (2,2) -- (0,4);
  \draw (0,0) -- (0,4);

  \draw (-2,2) -- (-2.5,2);
  \draw (0,0) -- (-0.5,-0.5);
  \draw (0,0) -- (0.5,-0.5);
  \fill (-0.25, -0.5) circle (0.02);
  \fill (0,-0.5) circle (0.02);
  \fill (0.25, -0.5) circle (0.02);

  \draw (0,4) -- (-0.5,4.5);
  \draw (0,4) -- (0.5,4.5);
  \fill (-0.25, 4.5) circle (0.02);
  \fill (0,4.5) circle (0.02);
  \fill (0.25, 4.5) circle (0.02);

  \draw[below right] (1,1) node {$a$};
  \draw[above right] (1,3) node {$a$};
  \draw[right] (0,2) node {$b$};
  \draw[above left] (-1,3) node {$c$};
  \draw[below left] (-1,1) node {$a$};
  \draw[left] (-2.5,2) node {$a$};

  \draw[very thick, dashed, ->] (3.5, 2.0) -- (4.5, 2.0);

  \begin{scope}[shift={(8.0, 0.0)}]
    \foreach \x/\y/\z in {0/0/$w$, 2/2/$u$, 0/4/$z$}
    {
      \fill (\x,\y) circle (0.1);
      \draw[right] (\x,\y) node {\z};
    }

    \fill (-2,2) circle (0.1);
    \draw[above] (-2,2) node {$v$};

    \draw (0,0) -- (2,2);
    \draw (0,4) -- (-2,2) -- (0,0);
    \draw (0,0) -- (0,4);

    \draw (-2,2) -- (-2.5,2);
    \draw (0,0) -- (-0.5,-0.5);
    \draw (0,0) -- (0.5,-0.5);
    \fill (-0.25, -0.5) circle (0.02);
    \fill (0,-0.5) circle (0.02);
    \fill (0.25, -0.5) circle (0.02);

    \draw (0,4) -- (-0.5,4.5);
    \draw (0,4) -- (0.5,4.5);
    \fill (-0.25, 4.5) circle (0.02);
    \fill (0,4.5) circle (0.02);
    \fill (0.25, 4.5) circle (0.02);
    \draw (2,2) -- (0,4);

    \draw[above right] (1,3) node {$b$};
    \draw[below right] (1,1) node {$a$};
    \draw[right] (0,2) node {$a$};
    \draw[above left] (-1,3) node {$c$};
    \draw[below left] (-1,1) node {$b$};
    \draw[left] (-2.5,2) node {$a$};
  \end{scope} 
\end{tikzpicture}
}}
\caption{Case $2.2$ of the proof of Lemma~\ref{lem:23triangles}.}
\end{center}
\end{figure}
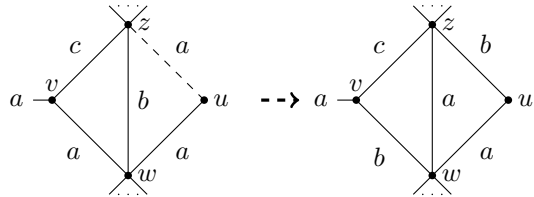

  $Case\ 2.3$ $C(wz)=b\not =a$, $C(vz)=a$, $C(vw)=c\not =a$.
  Just note that if we uncolor $wu$, color $zu$ with $a$ and swap the names of vertices $z$ and $w$ we arrive at Case 2.2.

  $Case\ 2.4$ $C(wz)=b\not =a$, $C(vz)=a$, $C(vw)=a$. 
  Since there is only one outside edge incident to the vertex $v$
  there can not be simultaneously $(vz,b)$-path and $(vw,b)$-path.
  Because of the symmetry we may assume w.l.o.g. that there is
  no $(vw,b)$-path. In this case we recolor $C(wz)=a$,$C(vw)=b$
  and put $C(uz)=b$.

\begin{figure}[h]
\begin{center}
{\small{
\begin{tikzpicture}[scale=0.5]

  \foreach \x/\y/\z in {0/0/$w$, 2/2/$u$, 0/4/$z$}
  {
    \fill (\x,\y) circle (0.1);
    \draw[right] (\x,\y) node {\z};
  }

  \fill (-2,2) circle (0.1);
  \draw[above] (-2,2) node {$v$};

  \draw (0,0) -- (2,2);
  \draw (0,4) -- (-2,2) -- (0,0);
  \draw[dashed] (2,2) -- (0,4);
  \draw (0,0) -- (0,4);

  \draw (-2,2) -- (-2.5,2);
  \draw (0,0) -- (-0.5,-0.5);
  \draw (0,0) -- (0.5,-0.5);
  \fill (-0.25, -0.5) circle (0.02);
  \fill (0,-0.5) circle (0.02);
  \fill (0.25, -0.5) circle (0.02);

  \draw (0,4) -- (-0.5,4.5);
  \draw (0,4) -- (0.5,4.5);
  \fill (-0.25, 4.5) circle (0.02);
  \fill (0,4.5) circle (0.02);
  \fill (0.25, 4.5) circle (0.02);

  \draw[below right] (1,1) node {$a$};
  \draw[above right] (1,3) node {$a$};
  \draw[right] (0,2) node {$b$};
  \draw[above left] (-1,3) node {$a$};
  \draw[below left] (-1,1) node {$a$};

  \draw[very thick, dashed, ->] (3.5, 2.0) -- (4.5, 2.0);

  \begin{scope}[shift={(8.0, 0.0)}]
    \foreach \x/\y/\z in {0/0/$w$, 2/2/$u$, 0/4/$z$}
    {
      \fill (\x,\y) circle (0.1);
      \draw[right] (\x,\y) node {\z};
    }

    \fill (-2,2) circle (0.1);
    \draw[above] (-2,2) node {$v$};

    \draw (0,0) -- (2,2);
    \draw (0,4) -- (-2,2) -- (0,0);
    \draw (0,0) -- (0,4);

    \draw (-2,2) -- (-2.5,2);
    \draw (0,0) -- (-0.5,-0.5);
    \draw (0,0) -- (0.5,-0.5);
    \fill (-0.25, -0.5) circle (0.02);
    \fill (0,-0.5) circle (0.02);
    \fill (0.25, -0.5) circle (0.02);

    \draw (0,4) -- (-0.5,4.5);
    \draw (0,4) -- (0.5,4.5);
    \fill (-0.25, 4.5) circle (0.02);
    \fill (0,4.5) circle (0.02);
    \fill (0.25, 4.5) circle (0.02);
    \draw (2,2) -- (0,4);

    \draw[above right] (1,3) node {$b$};
    \draw[below right] (1,1) node {$a$};
    \draw[right] (0,2) node {$a$};
    \draw[above left] (-1,3) node {$a$};
    \draw[below left] (-1,1) node {$b$};
  \end{scope} 
\end{tikzpicture}
}}
\caption{Case $2.4$ of the proof of Lemma~\ref{lem:23triangles}.}
\end{center}
\end{figure}
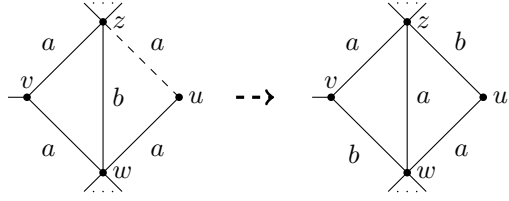

\end{proof}

\begin{figure}[ht]
	\centering
		\includegraphics{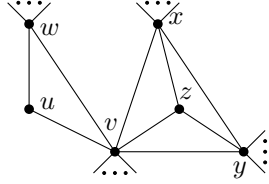}
		\caption{The configuration from Lemma~\ref{lem_2with3}.}
		\label{fig_2with3}
\end{figure}

\begin{lemma}
	\label{lem_2with3}
	$G$ does not contain the configuration in Fig.~\ref{fig_2with3}.
\end{lemma}

\begin{proof}
	Consider the configuration and the labeling of vertices as in Fig.~\ref{fig_2with3}. 
	By the minimality of $G$, there exists a $k$-linear coloring $C$ of $G - uv$. We show how to extend $C$ to $G$ and
	obtain a contradiction on the minimality. The only non-colored edge is $uv$. Let $a$ be the free color at $v$. We may assume that $a = C(uw)$
	and that there exists $(vu,a)$-path, for otherwise we can color $uv$ with $a$ without introducing a monochromatic cycle. 
	
	$Case\ 1.$ $C(vz)\ne a$. 
	We color $uv$ with $C(vz)$ and we uncolor $vz$, obtaining a $k$-linear coloring of $G-vz$ with a $(vu,a)$-path.
	We can assume that at least one of $zx$, $zy$ is not colored with $a$, for otherwise we just recolor both $zx$ and $zy$ to $C(xy)$ and
	$xy$ to $a$ and we obtain another $k$-linear coloring of $G-vz$ with none of $zx$, $zy$ colored with $a$. 
	Hence, we can color $vz$ with $a$ and we do not introduce a monochromatic cycle because $vz$ is on an $a$-path which ends at $u$.
	
	$Case\ 2.$ $C(vz)= a$. 
	Since there is a $(vu,a)$-path, $C(zx)=a$ or $C(zy)=a$. W.l.o.g.\ assume $C(zx)=a$. 
	Let $C(zy)=b$. We can assume that $C(xv)=b$ for otherwise we color both $zx$, $vz$ with $C(xv)$ and $xv$ with $a$ and we arrive at Case 1.
	Let $c=C(vw)$. If $c \neq b$, or $c = b$ and there is no $(xy, b)$-path, we can color both $uw$ and $vz$ with $c$, and both $vw$, $uv$ with $a$. 
	Hence, we assume that $C(vw) = b$ and there is an $(xy,b)$-path.
	Let $d=C(vy)$. Note that $a,b,d$ are pairwise distinct.
	Then we recolor the edges as in Fig.~\ref{lem_2with3gnote}, 
\begin{figure}[ht]
	\centering
		\includegraphics{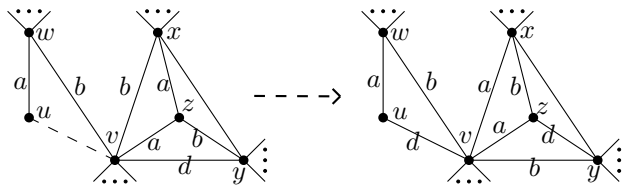}
		\caption{Proof of Lemma~\ref{lem_2with3}, Case 2.}
		\label{lem_2with3gnote}
\end{figure}	
	that is $vx$ with $a$, $vy$ and $xz$ with $b$, and $uv$, $yz$ with $d$. 
	We do not introduce any $d$-cycle, because $d\in C_1(u)$ and $d\in C_1(z)$.
	We do not introduce any $a$-cycle, because both $vx$ and $vz$ are on an $a$-path which ends at $z$.
	Finally, we do not introduce any $b$-cycle, because there is a $(xy,b)$-path so $vy$ is on a $b$-path which ends at $z$.
\end{proof}

From Lemmas~\ref{lem:two-pairs}, \ref{lem:23triangles} and \ref{lem_2with3} we immediately obtain the following corollary.

\begin{corollary}
\label{cor:atmostonetrianglewith3vertex}
If a vertex $v$ is incident to a triangle with a $2$-vertex then every $3$-vertex adjacent to $v$ is incident 
to at most one triangle.
\end{corollary}

\begin{lemma}
	\label{lem_cubicsmall}
	If $k\ge 3$ and $G$ contains a vertex $v$ of degree at most $2k-1$ with two $3$-neighbors then the neighbors of any $3$-neighbor of $v$ are pairwise nonadjacent.
\end{lemma}

\ignore{
\begin{figure}[ht]
\centering
	\includegraphics{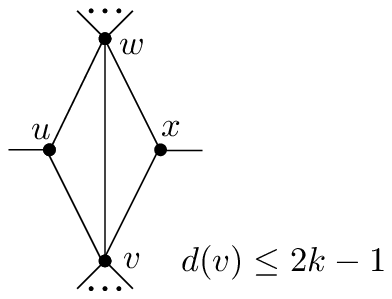}
	\caption{The configuration from Lemma~\ref{lem_cubicsmall}.}
	\label{fig_cubicsmall}
\end{figure}
}

\begin{proof}
    Let $x$ and $y$ be two $3$-neighbors of $v$ and assume that a pair of neighbors of $x$ or $y$, say of $x$, is adjacent.
	Let $x_1$, $x_2$ (resp. $y_1$, $y_2$) be the neighbors of $x$ (resp. $y$) distinct from $v$.
	By Lemma~\ref{lem_edge}, $y\not\in\{x_1,x_2\}$ and $x\not\in\{y_1,y_2\}$.
	
	Let $G'=G-vx$. Since $G$ is minimal, there exists a $k$-coloring $C$ of $G'$. In what follows, we show that $C$
	can be extended to $G$. Color the edges of $G$ with the colors assigned in $G'$. The only non-colored edge is $vx$.
	
	$Case\ 1.$ $C_0(v)\ne \emptyset$. Let $a$ be an element of $C_0(v)$. 
	We immediately infer that $C(xx_1) = C(xx_2) = a$, otherwise we color
	$vx$ with $a$ and introduce no monochromatic cycle, since no other edge incident to $v$ is colored by $a$. 
	If $x_1$ is adjacent with $v$, we color $xx_1$ and $vx$ with $C(x_1v)$ and $x_1v$ with $a$ and we obtain a $k$-linear coloring of $G$.
	We proceed similarly when $x_2$ is adjacent with $v$.
	Finally, when $x_1$ is adjacent with $x_2$, we just color $xx_1$ and $xx_2$ with $C(x_1x_2)$, and both $vx$ and $x_1x_2$ with $a$. 
	
    $Case\ 2.$ $C_0(v)= \emptyset$. Then $|C_1(v)|\ge 2$. Let $a,b$ be two distinct elements of $C_1(v)$. 	
	Observe that if $C(xx_1)$ and $C(xx_2)$ are both distinct from one of $a$ or $b$, we color $vx$ with that color. 
	Hence, we may assume that, without loss of generality, $C(xx_1) = a$ and $C(xx_2) = b$. 
	There exist also $(vx,a)$- and $(vx,b)$-paths, otherwise we color $vx$ with $a$ or $b$ without introducing a monochromatic cycle. 
	
	Let $c=C(vy)$. We color $vx$ with $c$. Next, we color $vy$ with $a$ if $a\not\in C_2(y)$ or with $b$ otherwise.
	It is easy to check that each color induces a graph of maximum degree 2. It suffices to check that neither $vx$ nor $vy$ belong
	to a monochromatic cycle.
	If $c=a$, then, since there is a $(vx,a)$-path, $C(yy_1)=a$ or $C(yy_2)=a$, so we colored $vy$ with $b$. 
	Hence, $vx$ is on an $a$-path ending at $v$ and since there is a $(vx,b)$-path, $vy$ is on a $b$-path ending at $x$.
	Now assume $c=b$. If we colored $vy$ with $b$ it means that $C(yy_1)=C(yy_2)=a$, so both $vx$ and $vy$ are on a $b$-path which ends at $y$.
	Otherwise, $vx$ is on a $b$-path ending at $v$ and since there is a $(vx,a)$-path, $vy$ is on a $a$-path ending at $x$.
	Finally, if $c\not\in\{a,b\}$, edge $vx$ is on a $c$-path ending at $x$ and edge $vy$ is on a monochromatic path ending at $x$.
\end{proof}

\begin{lemma}
  \label{lem:two_cubic}
	$G$ does not contain the configuration in Fig.~\ref{fig:two_cubic}.
\end{lemma}
\begin{figure}[hptb]
	\begin{center}
  {\small{
  \begin{tikzpicture}[scale=0.5]

      \foreach \x/\y/\z in {0/0, 0/4, -1.5/2, -4/2}
    {
      \fill (\x,\y) circle (0.1);
      \fill (-\x,\y) circle (0.1);
    }

    \foreach \x in {-1, 1}
    {
      \draw (\x*4,2) -- (0,0) -- (0,4) -- (\x * 1.5, 2) -- (0,0) -- (\x * 4, 2) -- (0,4);
    }
    \draw (4,2) -- (1.5,2);
    \draw (-4,2) -- (-1.5,2);

    \begin{scope}[shift={(0.0, 2.0)}]

      \foreach \x in {-1,1} 
      {
        \draw (0, \x * 2) -- (-0.5,\x * 2.5);
        \draw (0,\x * 2) -- (0.5, \x * 2.5);
        \fill (-0.25, \x * 2.5) circle (0.02);
        \fill (0,\x * 2.5) circle (0.02);
        \fill (0.25, \x * 2.5) circle (0.02);
      }

      \foreach \x in {-1,1} 
      {
        \draw (\x * 4, 0) -- (\x * 4.5, 0.5);
        \draw (\x * 4, 0) -- (\x * 4.5, -0.5);

        \fill (\x * 4.5, 0.25) circle (0.02);
        \fill (\x * 4.5, 0) circle (0.02);
        \fill (\x * 4.5, -0.25) circle (0.02);
      }
    \end{scope}

    \draw[above] (-4,2) node {$u$};
    \draw[above] (4,2) node {$y$};

    \draw[right] (0,4) node {$z$};
    \draw[right] (0,0) node {$w$};

    \draw[right] (-1.5,2) node {$v$};
    \draw[left] (1.5,2) node {$x$};

  \end{tikzpicture}
  }}
	\end{center}
	\caption{The configuration from Lemma~\ref{lem:two_cubic}.}
	\label{fig:two_cubic}
\end{figure}

\begin{proof}
Consider the configuration and the labeling of vertices as in Fig.~\ref{fig:two_cubic}.
Since $G$ is minimal, there exists a $k$-linear coloring $C$ of $G'=G-uv$.
Let $a$ be a free color at $u$. 
Now let us take into consideration edges $vz$ and $vw$.

$Case\ 1$ $C(vz)\not = a,C(vw)\not = a$. In this case we simply put $C(uv)=a$ and we are done.

$Case\ 2$ $C(vz)=a,C(vw)=a$. Let $b=C(zw)$. Obviously $b\not = a$ so we can recolor $C(vz)=C(vw)=b$,
$C(zw)=a$ and put $C(uv)=a$.

$Case\ 3$ Exactly one edge from the set $\{vz, vw\}$ has color $a$. Because of the symmetry we may
assume that $C(vw)=a$ and $C(vz)=b\not=a$. 
We additionally assume that there is a $(uw,a)$-path since otherwise we can put $C(uv)=a$ without
introducing a monochromatic cycle. 
Le us consider the edge $wx$.

$Case\ 3.1$ $C(vw)=a,C(vz)=b\not=a,C(wx)=c\not=a$ and there exists a $(uw,a)$-path (see Figure~\ref{fig:two_cubic_3_1}). 
We would like to swap colors on edges $wv$ and $wx$ thus in order to do so we consider subcases regarding the number of $a$-edges
incident to the vertex $x$.

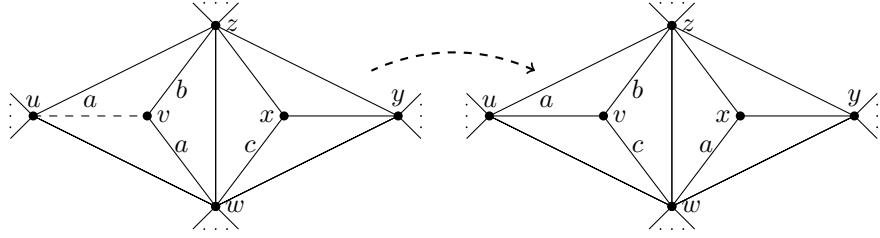
\begin{figure}[hptb]
	\begin{center}
  {\small{
  \begin{tikzpicture}[scale=0.6]

        \foreach \x/\y/\z in {0/0, 0/4, -1.5/2, -4/2}
    {
      \fill (\x,\y) circle (0.1);
      \fill (-\x,\y) circle (0.1);
    }

    \foreach \x in {-1, 1}
    {
      \draw (\x*4,2) -- (0,0) -- (0,4) -- (\x * 1.5, 2) -- (0,0) -- (\x * 4, 2) -- (0,4);
    }
    \draw (4,2) -- (1.5,2);
    \draw[dashed] (-4,2) -- (-1.5,2);

    \begin{scope}[shift={(0.0, 2.0)}]

      \foreach \x in {-1,1} 
      {
        \draw (0, \x * 2) -- (-0.5,\x * 2.5);
        \draw (0,\x * 2) -- (0.5, \x * 2.5);
        \fill (-0.25, \x * 2.5) circle (0.02);
        \fill (0,\x * 2.5) circle (0.02);
        \fill (0.25, \x * 2.5) circle (0.02);
      }

      \foreach \x in {-1,1} 
      {
        \draw (\x * 4, 0) -- (\x * 4.5, 0.5);
        \draw (\x * 4, 0) -- (\x * 4.5, -0.5);

        \fill (\x * 4.5, 0.25) circle (0.02);
        \fill (\x * 4.5, 0) circle (0.02);
        \fill (\x * 4.5, -0.25) circle (0.02);
      }
    \end{scope}

    \draw[above] (-4,2) node {$u$};
    \draw[above] (4,2) node {$y$};

    \draw[right] (0,4) node {$z$};
    \draw[right] (0,0) node {$w$};

    \draw[right] (-1.5,2) node {$v$};
    \draw[left] (1.5,2) node {$x$};

    \draw[above] (-2.75,2) node {$a$};
    \draw[above] (-0.75,1) node {$a$};
    \draw[below] (-0.75,3) node {$b$};
    \draw[above] (0.75,1) node {$c$};

  \end{tikzpicture}
  }}
	\end{center}
	\caption{Case $3.1$ in the proof of Lemma~\ref{lem:two_cubic}}
	\label{fig:two_cubic_3_1}
\end{figure}

$Case\ 3.1.1$ $C(vw)=a,C(vz)=b\not=a,C(wx)=c\not=a$, there exists a $(uw,a)$-path and there is at most one $a$-edge
incident to the vertex $x$. Let us swap colors of edges $vw$ and $wx$ as in Figure~\ref{fig:two_cubic_3_1_2}.
We know that the $a$-path that starts in the vertex $u$ reaches the vertex $w$ which means that it ends in the vertex $v$.
Since we have an assumption that there is at most one $a$-edge incident to the vertex $x$ it can not happen that
the $(uw,a)$-path goes through the vertex $x$ thus even if we connect those paths by swapping colors of edges $vw$ and $wx$
we do not introduce an $a$-cycle. We can only introduce a monochromatic cycle when $b=c$ and there exists a $(zw,b)$-path which does not
go through the edge $wx$.
We assume that this is the case since otherwise we are done. Let us consider the edge $xz$.

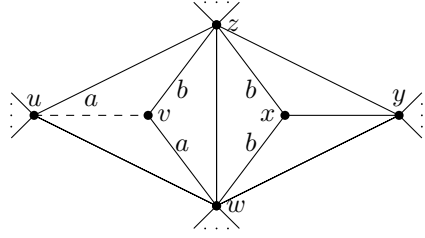
\begin{figure}[hptb]
	\begin{center}
  {\small{
  \begin{tikzpicture}[scale=0.6]

        \foreach \x/\y/\z in {0/0, 0/4, -1.5/2, -4/2}
    {
      \fill (\x,\y) circle (0.1);
      \fill (-\x,\y) circle (0.1);
    }

    \foreach \x in {-1, 1}
    {
      \draw (\x*4,2) -- (0,0) -- (0,4) -- (\x * 1.5, 2) -- (0,0) -- (\x * 4, 2) -- (0,4);
    }
    \draw (4,2) -- (1.5,2);
    \draw[dashed] (-4,2) -- (-1.5,2);

    \begin{scope}[shift={(0.0, 2.0)}]

      \foreach \x in {-1,1} 
      {
        \draw (0, \x * 2) -- (-0.5,\x * 2.5);
        \draw (0,\x * 2) -- (0.5, \x * 2.5);
        \fill (-0.25, \x * 2.5) circle (0.02);
        \fill (0,\x * 2.5) circle (0.02);
        \fill (0.25, \x * 2.5) circle (0.02);
      }

      \foreach \x in {-1,1} 
      {
        \draw (\x * 4, 0) -- (\x * 4.5, 0.5);
        \draw (\x * 4, 0) -- (\x * 4.5, -0.5);

        \fill (\x * 4.5, 0.25) circle (0.02);
        \fill (\x * 4.5, 0) circle (0.02);
        \fill (\x * 4.5, -0.25) circle (0.02);
      }
    \end{scope}

    \draw[above] (-4,2) node {$u$};
    \draw[above] (4,2) node {$y$};

    \draw[right] (0,4) node {$z$};
    \draw[right] (0,0) node {$w$};

    \draw[right] (-1.5,2) node {$v$};
    \draw[left] (1.5,2) node {$x$};

    \draw[above] (-2.75,2) node {$a$};
    \draw[above] (-0.75,1) node {$a$};
    \draw[below] (-0.75,3) node {$b$};
    \draw[above] (0.75,1) node {$c$};

    \draw[dashed,thick,<-] (7,3) arc (65:115:3*1.414);

    \begin{scope}[shift={(10.0, 0.0)}]

          \foreach \x/\y/\z in {0/0, 0/4, -1.5/2, -4/2}
    {
      \fill (\x,\y) circle (0.1);
      \fill (-\x,\y) circle (0.1);
    }

    \foreach \x in {-1, 1}
    {
      \draw (\x*4,2) -- (0,0) -- (0,4) -- (\x * 1.5, 2) -- (0,0) -- (\x * 4, 2) -- (0,4);
    }
    \draw (4,2) -- (1.5,2);
    \draw (-4,2) -- (-1.5,2);

    \begin{scope}[shift={(0.0, 2.0)}]

      \foreach \x in {-1,1} 
      {
        \draw (0, \x * 2) -- (-0.5,\x * 2.5);
        \draw (0,\x * 2) -- (0.5, \x * 2.5);
        \fill (-0.25, \x * 2.5) circle (0.02);
        \fill (0,\x * 2.5) circle (0.02);
        \fill (0.25, \x * 2.5) circle (0.02);
      }

      \foreach \x in {-1,1} 
      {
        \draw (\x * 4, 0) -- (\x * 4.5, 0.5);
        \draw (\x * 4, 0) -- (\x * 4.5, -0.5);

        \fill (\x * 4.5, 0.25) circle (0.02);
        \fill (\x * 4.5, 0) circle (0.02);
        \fill (\x * 4.5, -0.25) circle (0.02);
      }
    \end{scope}

    \draw[above] (-4,2) node {$u$};
    \draw[above] (4,2) node {$y$};

    \draw[right] (0,4) node {$z$};
    \draw[right] (0,0) node {$w$};

    \draw[right] (-1.5,2) node {$v$};
    \draw[left] (1.5,2) node {$x$};

      \draw[above] (-2.75,2) node {$a$};
      \draw[above] (-0.75,1) node {$c$};
      \draw[below] (-0.75,3) node {$b$};
      \draw[above] (0.75,1) node {$a$};

    \end{scope}

  \end{tikzpicture}
  }}
	\end{center}
	\caption{Case $3.1.1$ in the proof of Lemma~\ref{lem:two_cubic}}
	\label{fig:two_cubic_3_1_2}
\end{figure}

$Case\ 3.1.1.1$ $C(vw)=a,C(vz)=b\not=a,C(wx)=b,C(xz)=b$, there is at most one $a$-edge incident to the vertex $x$,
there exists a $(uw,a)$-path and there exists a $(zw,b)$-path which does not go through the edge $wx$.
Since the vertex $z$ has already two incident $b$-edges the last condition can not be satisfied, contradiction 
(see Figure~\ref{fig:two_cubic_3_1_2_1}).

\begin{figure}[hptb]
	\begin{center}
  {\small{
  \begin{tikzpicture}[scale=0.6]

        \foreach \x/\y/\z in {0/0, 0/4, -1.5/2, -4/2}
    {
      \fill (\x,\y) circle (0.1);
      \fill (-\x,\y) circle (0.1);
    }

    \foreach \x in {-1, 1}
    {
      \draw (\x*4,2) -- (0,0) -- (0,4) -- (\x * 1.5, 2) -- (0,0) -- (\x * 4, 2) -- (0,4);
    }
    \draw (4,2) -- (1.5,2);
    \draw[dashed] (-4,2) -- (-1.5,2);

    \begin{scope}[shift={(0.0, 2.0)}]

      \foreach \x in {-1,1} 
      {
        \draw (0, \x * 2) -- (-0.5,\x * 2.5);
        \draw (0,\x * 2) -- (0.5, \x * 2.5);
        \fill (-0.25, \x * 2.5) circle (0.02);
        \fill (0,\x * 2.5) circle (0.02);
        \fill (0.25, \x * 2.5) circle (0.02);
      }

      \foreach \x in {-1,1} 
      {
        \draw (\x * 4, 0) -- (\x * 4.5, 0.5);
        \draw (\x * 4, 0) -- (\x * 4.5, -0.5);

        \fill (\x * 4.5, 0.25) circle (0.02);
        \fill (\x * 4.5, 0) circle (0.02);
        \fill (\x * 4.5, -0.25) circle (0.02);
      }
    \end{scope}

    \draw[above] (-4,2) node {$u$};
    \draw[above] (4,2) node {$y$};

    \draw[right] (0,4) node {$z$};
    \draw[right] (0,0) node {$w$};

    \draw[right] (-1.5,2) node {$v$};
    \draw[left] (1.5,2) node {$x$};

    \draw[above] (-2.75,2) node {$a$};
    \draw[above] (-0.75,1) node {$a$};
    \draw[below] (-0.75,3) node {$b$};
    \draw[above] (0.75,1) node {$b$};
    \draw[below] (0.75,3) node {$b$};

  \end{tikzpicture}
  }}
	\end{center}
	\caption{Case $3.1.1.1$ in the proof of Lemma~\ref{lem:two_cubic}}
	\label{fig:two_cubic_3_1_2_1}
\end{figure}

$Case\ 3.1.1.2$ $C(vw)=a,C(vz)=b\not=a,C(wx)=b,C(xz)=c\not=b$, there is at most one $a$-edge incident to the vertex $x$,
there exists a $(uw,a)$-path and there exists a $(zw,b)$-path which does not go through the edge $wx$.
In this case we swap colors of two pairs of edges $\{zv,zx\}$ and $\{wv,wx\}$ as in Figure~\ref{fig:two_cubic_3_1_2_2}.
We show that no monochromatic cycle is introduced. 
Since there is a $(zw,b)$-path that does not go through $wx$, edges $zx$ and $vw$ are on the same $b$-path, which ends at $v$.
If $c=a$ and there is an $a$-cycle, it means that $C(xy)=a$ and there is a $(yz,a)$-path which does not go through $x$ --- but then
there is an $a$-cycle in $G'$, a contradiction.
If $c\ne a$, $vz$ is on a $c$-path that ends at $v$ and both $uv$ and $wx$ are on the same $a$-path, which ends at $v$.

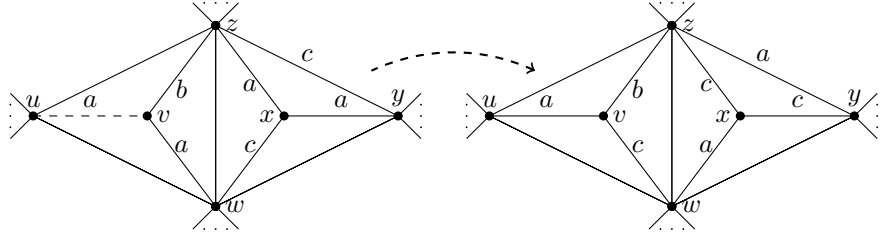
\begin{figure}[hptb]
	\begin{center}
  {\small{
  \begin{tikzpicture}[scale=0.6]

        \foreach \x/\y/\z in {0/0, 0/4, -1.5/2, -4/2}
    {
      \fill (\x,\y) circle (0.1);
      \fill (-\x,\y) circle (0.1);
    }

    \foreach \x in {-1, 1}
    {
      \draw (\x*4,2) -- (0,0) -- (0,4) -- (\x * 1.5, 2) -- (0,0) -- (\x * 4, 2) -- (0,4);
    }
    \draw (4,2) -- (1.5,2);
    \draw[dashed] (-4,2) -- (-1.5,2);

    \begin{scope}[shift={(0.0, 2.0)}]

      \foreach \x in {-1,1} 
      {
        \draw (0, \x * 2) -- (-0.5,\x * 2.5);
        \draw (0,\x * 2) -- (0.5, \x * 2.5);
        \fill (-0.25, \x * 2.5) circle (0.02);
        \fill (0,\x * 2.5) circle (0.02);
        \fill (0.25, \x * 2.5) circle (0.02);
      }

      \foreach \x in {-1,1} 
      {
        \draw (\x * 4, 0) -- (\x * 4.5, 0.5);
        \draw (\x * 4, 0) -- (\x * 4.5, -0.5);

        \fill (\x * 4.5, 0.25) circle (0.02);
        \fill (\x * 4.5, 0) circle (0.02);
        \fill (\x * 4.5, -0.25) circle (0.02);
      }
    \end{scope}

    \draw[above] (-4,2) node {$u$};
    \draw[above] (4,2) node {$y$};

    \draw[right] (0,4) node {$z$};
    \draw[right] (0,0) node {$w$};

    \draw[right] (-1.5,2) node {$v$};
    \draw[left] (1.5,2) node {$x$};

    \draw[above] (-2.75,2) node {$a$};
    \draw[above] (-0.75,1) node {$a$};
    \draw[below] (-0.75,3) node {$b$};
    \draw[above] (0.75,1) node {$b$};
    \draw[below] (0.75,3) node {$c$};

    \draw[dashed,thick,<-] (7,3) arc (65:115:3*1.414);

    \begin{scope}[shift={(10.0, 0.0)}]

          \foreach \x/\y/\z in {0/0, 0/4, -1.5/2, -4/2}
    {
      \fill (\x,\y) circle (0.1);
      \fill (-\x,\y) circle (0.1);
    }

    \foreach \x in {-1, 1}
    {
      \draw (\x*4,2) -- (0,0) -- (0,4) -- (\x * 1.5, 2) -- (0,0) -- (\x * 4, 2) -- (0,4);
    }
    \draw (4,2) -- (1.5,2);
    \draw (-4,2) -- (-1.5,2);

    \begin{scope}[shift={(0.0, 2.0)}]

      \foreach \x in {-1,1} 
      {
        \draw (0, \x * 2) -- (-0.5,\x * 2.5);
        \draw (0,\x * 2) -- (0.5, \x * 2.5);
        \fill (-0.25, \x * 2.5) circle (0.02);
        \fill (0,\x * 2.5) circle (0.02);
        \fill (0.25, \x * 2.5) circle (0.02);
      }

      \foreach \x in {-1,1} 
      {
        \draw (\x * 4, 0) -- (\x * 4.5, 0.5);
        \draw (\x * 4, 0) -- (\x * 4.5, -0.5);

        \fill (\x * 4.5, 0.25) circle (0.02);
        \fill (\x * 4.5, 0) circle (0.02);
        \fill (\x * 4.5, -0.25) circle (0.02);
      }
    \end{scope}

    \draw[above] (-4,2) node {$u$};
    \draw[above] (4,2) node {$y$};

    \draw[right] (0,4) node {$z$};
    \draw[right] (0,0) node {$w$};

    \draw[right] (-1.5,2) node {$v$};
    \draw[left] (1.5,2) node {$x$};

      \draw[above] (-2.75,2) node {$a$};
      \draw[above] (-0.75,1) node {$b$};
      \draw[below] (-0.75,3) node {$c$};
      \draw[above] (0.75,1) node {$a$};
      \draw[below] (0.75,3) node {$b$};

    \end{scope}

  \end{tikzpicture}
  }}
	\end{center}
	\caption{Case $3.1.1.2$ in the proof of Lemma~\ref{lem:two_cubic}}
	\label{fig:two_cubic_3_1_2_2}
\end{figure}

$Case\ 3.1.2$ $C(vw)=a,C(vz)=b\not=a,C(wx)=c\not=a,C(zx)=C(xy)=a$ and there exists a $(uw,a)$-path.
We may assume that $C(zy)=c$ since if $C(zy)=d\not=c$ we can recolor $C(zx)=C(xy)=d$ and $C(zy)=a$
which would result in the same situation as in the already solved Case $3.1.1$.
Now we try to recolor $C(vw)=c,C(wx)=a,C(zy)=a,C(zx)=C(xy)=c$ as in Figure~\ref{fig:two_cubic_3_1_1}. 
Now let us consider cases:

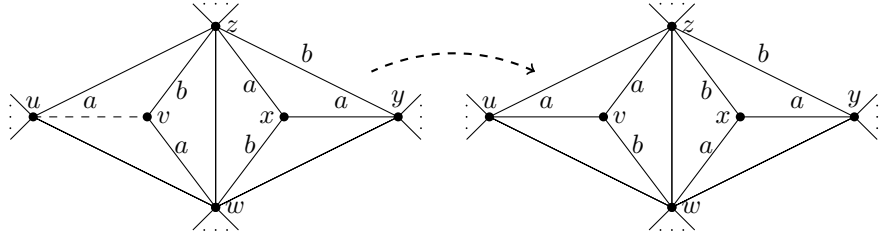
\begin{figure}[hptb]
	\begin{center}
  {\small{
  \begin{tikzpicture}[scale=0.6]

        \foreach \x/\y/\z in {0/0, 0/4, -1.5/2, -4/2}
    {
      \fill (\x,\y) circle (0.1);
      \fill (-\x,\y) circle (0.1);
    }

    \foreach \x in {-1, 1}
    {
      \draw (\x*4,2) -- (0,0) -- (0,4) -- (\x * 1.5, 2) -- (0,0) -- (\x * 4, 2) -- (0,4);
    }
    \draw (4,2) -- (1.5,2);
    \draw[dashed] (-4,2) -- (-1.5,2);

    \begin{scope}[shift={(0.0, 2.0)}]

      \foreach \x in {-1,1} 
      {
        \draw (0, \x * 2) -- (-0.5,\x * 2.5);
        \draw (0,\x * 2) -- (0.5, \x * 2.5);
        \fill (-0.25, \x * 2.5) circle (0.02);
        \fill (0,\x * 2.5) circle (0.02);
        \fill (0.25, \x * 2.5) circle (0.02);
      }

      \foreach \x in {-1,1} 
      {
        \draw (\x * 4, 0) -- (\x * 4.5, 0.5);
        \draw (\x * 4, 0) -- (\x * 4.5, -0.5);

        \fill (\x * 4.5, 0.25) circle (0.02);
        \fill (\x * 4.5, 0) circle (0.02);
        \fill (\x * 4.5, -0.25) circle (0.02);
      }
    \end{scope}

    \draw[above] (-4,2) node {$u$};
    \draw[above] (4,2) node {$y$};

    \draw[right] (0,4) node {$z$};
    \draw[right] (0,0) node {$w$};

    \draw[right] (-1.5,2) node {$v$};
    \draw[left] (1.5,2) node {$x$};

    \draw[above] (-2.75,2) node {$a$};
    \draw[above] (-0.75,1) node {$a$};
    \draw[below] (-0.75,3) node {$b$};
    \draw[above] (0.75,1) node {$c$};

    \draw[above] (2.75,2) node {$a$};
    \draw[below] (0.75,3) node {$a$};
    \draw[above] (2,3) node {$c$};

    \draw[dashed,thick,<-] (7,3) arc (65:115:3*1.414);

    \begin{scope}[shift={(10.0, 0.0)}]
          \foreach \x/\y/\z in {0/0, 0/4, -1.5/2, -4/2}
    {
      \fill (\x,\y) circle (0.1);
      \fill (-\x,\y) circle (0.1);
    }

    \foreach \x in {-1, 1}
    {
      \draw (\x*4,2) -- (0,0) -- (0,4) -- (\x * 1.5, 2) -- (0,0) -- (\x * 4, 2) -- (0,4);
    }
    \draw (4,2) -- (1.5,2);
    \draw (-4,2) -- (-1.5,2);

    \begin{scope}[shift={(0.0, 2.0)}]

      \foreach \x in {-1,1} 
      {
        \draw (0, \x * 2) -- (-0.5,\x * 2.5);
        \draw (0,\x * 2) -- (0.5, \x * 2.5);
        \fill (-0.25, \x * 2.5) circle (0.02);
        \fill (0,\x * 2.5) circle (0.02);
        \fill (0.25, \x * 2.5) circle (0.02);
      }

      \foreach \x in {-1,1} 
      {
        \draw (\x * 4, 0) -- (\x * 4.5, 0.5);
        \draw (\x * 4, 0) -- (\x * 4.5, -0.5);

        \fill (\x * 4.5, 0.25) circle (0.02);
        \fill (\x * 4.5, 0) circle (0.02);
        \fill (\x * 4.5, -0.25) circle (0.02);
      }
    \end{scope}

    \draw[above] (-4,2) node {$u$};
    \draw[above] (4,2) node {$y$};

    \draw[right] (0,4) node {$z$};
    \draw[right] (0,0) node {$w$};

    \draw[right] (-1.5,2) node {$v$};
    \draw[left] (1.5,2) node {$x$};

      \draw[above] (-2.75,2) node {$a$};
      \draw[above] (-0.75,1) node {$c$};
      \draw[below] (-0.75,3) node {$b$};
      \draw[above] (0.75,1) node {$a$};

      \draw[above] (2.75,2) node {$c$};
      \draw[below] (0.75,3) node {$c$};
      \draw[above] (2,3) node {$a$};
    \end{scope}

  \end{tikzpicture}
  }}
	\end{center}
	\caption{Case $3.1.2$ in the proof of Lemma~\ref{lem:two_cubic}}
	\label{fig:two_cubic_3_1_1}
\end{figure}

$Case\ 3.1.2.1$ We do not introduce a monochromatic cycle and we are done.

$Case\ 3.1.2.2$ We introduce a monochromatic cycle which means that $b=c$ and there exists a $(wy,b)$-path, hence
the set of assumptions is $C(vw)=a,C(vz)=b\not=a,C(wx)=b,C(zx)=C(xy)=a,C(zy)=b$, there exists a $(uw,a)$-path and
there exists a $(wy,b)$-path. 
In this case we consider the $(uw,a)$-path and branch on the relation between this path and vertices $z$, $x$, and $y$.

$Case\ 3.1.2.2.1$ $C(vw)=a,C(vz)=b\not=a,C(wx)=b,C(zx)=C(xy)=a,C(zy)=b$, there exists a $(uw,a)$-path,
there exists a $(wy,b)$-path and the $(uw,a)$-path does not go through the vertex $x$.
In this case we swap colors on two pairs of edges $\{zv,zx\}$ and $\{wv,wx\}$ as in Figure~\ref{fig:two_cubic_3_1_1_2_1}
and in this way we do not introduce any monochromatic cycle since we join distinct $a$-paths.

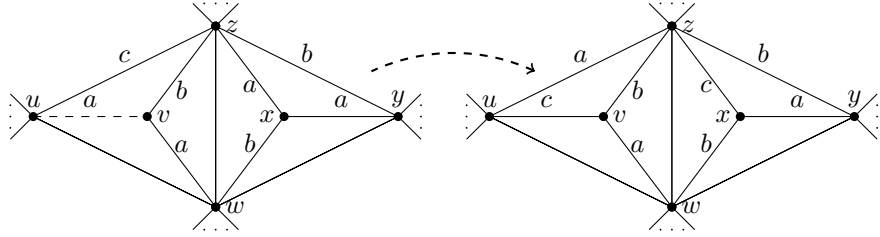
\begin{figure}[hptb]
	\begin{center}
  {\small{
  \begin{tikzpicture}[scale=0.6]

        \foreach \x/\y/\z in {0/0, 0/4, -1.5/2, -4/2}
    {
      \fill (\x,\y) circle (0.1);
      \fill (-\x,\y) circle (0.1);
    }

    \foreach \x in {-1, 1}
    {
      \draw (\x*4,2) -- (0,0) -- (0,4) -- (\x * 1.5, 2) -- (0,0) -- (\x * 4, 2) -- (0,4);
    }
    \draw (4,2) -- (1.5,2);
    \draw[dashed] (-4,2) -- (-1.5,2);

    \begin{scope}[shift={(0.0, 2.0)}]

      \foreach \x in {-1,1} 
      {
        \draw (0, \x * 2) -- (-0.5,\x * 2.5);
        \draw (0,\x * 2) -- (0.5, \x * 2.5);
        \fill (-0.25, \x * 2.5) circle (0.02);
        \fill (0,\x * 2.5) circle (0.02);
        \fill (0.25, \x * 2.5) circle (0.02);
      }

      \foreach \x in {-1,1} 
      {
        \draw (\x * 4, 0) -- (\x * 4.5, 0.5);
        \draw (\x * 4, 0) -- (\x * 4.5, -0.5);

        \fill (\x * 4.5, 0.25) circle (0.02);
        \fill (\x * 4.5, 0) circle (0.02);
        \fill (\x * 4.5, -0.25) circle (0.02);
      }
    \end{scope}

    \draw[above] (-4,2) node {$u$};
    \draw[above] (4,2) node {$y$};

    \draw[right] (0,4) node {$z$};
    \draw[right] (0,0) node {$w$};

    \draw[right] (-1.5,2) node {$v$};
    \draw[left] (1.5,2) node {$x$};

    \draw[above] (-2.75,2) node {$a$};
    \draw[above] (-0.75,1) node {$a$};
    \draw[below] (-0.75,3) node {$b$};
    \draw[above] (0.75,1) node {$b$};

    \draw[above] (2.75,2) node {$a$};
    \draw[below] (0.75,3) node {$a$};
    \draw[above] (2,3) node {$b$};

    \draw[dashed,thick,<-] (7,3) arc (65:115:3*1.414);

    \begin{scope}[shift={(10.0, 0.0)}]
          \foreach \x/\y/\z in {0/0, 0/4, -1.5/2, -4/2}
    {
      \fill (\x,\y) circle (0.1);
      \fill (-\x,\y) circle (0.1);
    }

    \foreach \x in {-1, 1}
    {
      \draw (\x*4,2) -- (0,0) -- (0,4) -- (\x * 1.5, 2) -- (0,0) -- (\x * 4, 2) -- (0,4);
    }
    \draw (4,2) -- (1.5,2);
    \draw (-4,2) -- (-1.5,2);

    \begin{scope}[shift={(0.0, 2.0)}]

      \foreach \x in {-1,1} 
      {
        \draw (0, \x * 2) -- (-0.5,\x * 2.5);
        \draw (0,\x * 2) -- (0.5, \x * 2.5);
        \fill (-0.25, \x * 2.5) circle (0.02);
        \fill (0,\x * 2.5) circle (0.02);
        \fill (0.25, \x * 2.5) circle (0.02);
      }

      \foreach \x in {-1,1} 
      {
        \draw (\x * 4, 0) -- (\x * 4.5, 0.5);
        \draw (\x * 4, 0) -- (\x * 4.5, -0.5);

        \fill (\x * 4.5, 0.25) circle (0.02);
        \fill (\x * 4.5, 0) circle (0.02);
        \fill (\x * 4.5, -0.25) circle (0.02);
      }
    \end{scope}

    \draw[above] (-4,2) node {$u$};
    \draw[above] (4,2) node {$y$};

    \draw[right] (0,4) node {$z$};
    \draw[right] (0,0) node {$w$};

    \draw[right] (-1.5,2) node {$v$};
    \draw[left] (1.5,2) node {$x$};

      \draw[above] (-2.75,2) node {$a$};
      \draw[above] (-0.75,1) node {$b$};
      \draw[below] (-0.75,3) node {$a$};
      \draw[above] (0.75,1) node {$a$};

      \draw[above] (2.75,2) node {$a$};
      \draw[below] (0.75,3) node {$b$};
      \draw[above] (2,3) node {$b$};
    \end{scope}

  \end{tikzpicture}
  }}
	\end{center}
	\caption{Case $3.1.2.2.1$ in the proof of Lemma~\ref{lem:two_cubic}}
	\label{fig:two_cubic_3_1_1_2_1}
\end{figure}

$Case\ 3.1.2.2.2$ $C(vw)=a,C(vz)=b\not=a,C(wx)=b,C(zx)=C(xy)=a,C(zy)=b$, there exists 
a $(wy,b)$-path, there exists a $(wz,a)$-path which does not go through the vertex $x$
and there exists a $(yu,a)$-path which does not go through the vertex $x$.
Let $c=C(uz)$. Since the vertex $z$ already has two incident $b$-edges we have $c\not=b$.
Moreover the $(yu,a)$-path does not go through the vertex $x$ thus $c\not=a$.
In this case we can recolor $C(uz)=a,C(uv)=c,C(zx)=c$ as in Figure~\ref{fig:two_cubic_3_1_1_2_2}
and we do not introduce any monochromatic cycle because colors $a,b,c$ are pairwise distinct
and all the $a$-edges of the configuration are on an $a$-path which ends at $v$ and $x$.

\begin{figure}[hptb]
	\begin{center}
  {\small{
  \begin{tikzpicture}[scale=0.6]

        \foreach \x/\y/\z in {0/0, 0/4, -1.5/2, -4/2}
    {
      \fill (\x,\y) circle (0.1);
      \fill (-\x,\y) circle (0.1);
    }

    \foreach \x in {-1, 1}
    {
      \draw (\x*4,2) -- (0,0) -- (0,4) -- (\x * 1.5, 2) -- (0,0) -- (\x * 4, 2) -- (0,4);
    }
    \draw (4,2) -- (1.5,2);
    \draw[dashed] (-4,2) -- (-1.5,2);

    \begin{scope}[shift={(0.0, 2.0)}]

      \foreach \x in {-1,1} 
      {
        \draw (0, \x * 2) -- (-0.5,\x * 2.5);
        \draw (0,\x * 2) -- (0.5, \x * 2.5);
        \fill (-0.25, \x * 2.5) circle (0.02);
        \fill (0,\x * 2.5) circle (0.02);
        \fill (0.25, \x * 2.5) circle (0.02);
      }

      \foreach \x in {-1,1} 
      {
        \draw (\x * 4, 0) -- (\x * 4.5, 0.5);
        \draw (\x * 4, 0) -- (\x * 4.5, -0.5);

        \fill (\x * 4.5, 0.25) circle (0.02);
        \fill (\x * 4.5, 0) circle (0.02);
        \fill (\x * 4.5, -0.25) circle (0.02);
      }
    \end{scope}

    \draw[above] (-4,2) node {$u$};
    \draw[above] (4,2) node {$y$};

    \draw[right] (0,4) node {$z$};
    \draw[right] (0,0) node {$w$};

    \draw[right] (-1.5,2) node {$v$};
    \draw[left] (1.5,2) node {$x$};

    \draw[above] (-2.75,2) node {$a$};
    \draw[above] (-0.75,1) node {$a$};
    \draw[below] (-0.75,3) node {$b$};
    \draw[above] (0.75,1) node {$b$};

    \draw[above] (2.75,2) node {$a$};
    \draw[below] (0.75,3) node {$a$};
    \draw[above] (2,3) node {$b$};
    \draw[above] (-2,3) node {$c$};

    \draw[dashed,thick,<-] (7,3) arc (65:115:3*1.414);

    \begin{scope}[shift={(10.0, 0.0)}]
          \foreach \x/\y/\z in {0/0, 0/4, -1.5/2, -4/2}
    {
      \fill (\x,\y) circle (0.1);
      \fill (-\x,\y) circle (0.1);
    }

    \foreach \x in {-1, 1}
    {
      \draw (\x*4,2) -- (0,0) -- (0,4) -- (\x * 1.5, 2) -- (0,0) -- (\x * 4, 2) -- (0,4);
    }
    \draw (4,2) -- (1.5,2);
    \draw (-4,2) -- (-1.5,2);

    \begin{scope}[shift={(0.0, 2.0)}]

      \foreach \x in {-1,1} 
      {
        \draw (0, \x * 2) -- (-0.5,\x * 2.5);
        \draw (0,\x * 2) -- (0.5, \x * 2.5);
        \fill (-0.25, \x * 2.5) circle (0.02);
        \fill (0,\x * 2.5) circle (0.02);
        \fill (0.25, \x * 2.5) circle (0.02);
      }

      \foreach \x in {-1,1} 
      {
        \draw (\x * 4, 0) -- (\x * 4.5, 0.5);
        \draw (\x * 4, 0) -- (\x * 4.5, -0.5);

        \fill (\x * 4.5, 0.25) circle (0.02);
        \fill (\x * 4.5, 0) circle (0.02);
        \fill (\x * 4.5, -0.25) circle (0.02);
      }
    \end{scope}

    \draw[above] (-4,2) node {$u$};
    \draw[above] (4,2) node {$y$};

    \draw[right] (0,4) node {$z$};
    \draw[right] (0,0) node {$w$};

    \draw[right] (-1.5,2) node {$v$};
    \draw[left] (1.5,2) node {$x$};

      \draw[above] (-2.75,2) node {$c$};
      \draw[above] (-0.75,1) node {$a$};
      \draw[below] (-0.75,3) node {$b$};
      \draw[above] (0.75,1) node {$b$};

      \draw[above] (2.75,2) node {$a$};
      \draw[below] (0.75,3) node {$c$};
      \draw[above] (2,3) node {$b$};
      \draw[above] (-2,3) node {$a$};
    \end{scope}

  \end{tikzpicture}
  }}
	\end{center}
	\caption{Case $3.1.2.2.2$ in the proof of Lemma~\ref{lem:two_cubic}}
	\label{fig:two_cubic_3_1_1_2_2}
\end{figure}
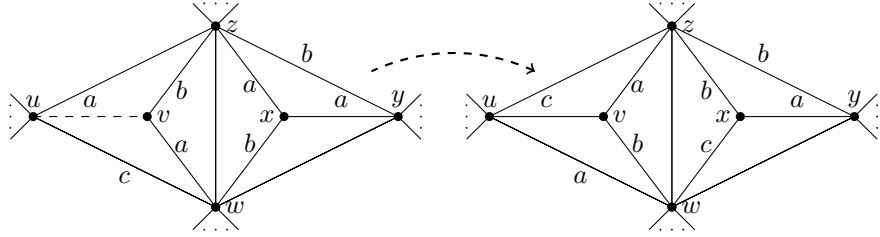

$Case\ 3.1.2.2.3$ $C(vw)=a,C(vz)=b\not=a,C(wx)=b,C(zx)=C(xy)=a,C(zy)=b$, there exists 
a $(wy,b)$-path, there exists a $(wy,a)$-path which does not go through the vertex $x$
and there exists a $(zu,a)$-path which does not go through the vertex $x$.
Let $c=C(uw)$. Because of the $(zu,a)$-path we have $c\not=a$ thus we branch
into cases were $c\ne b$ and $c=b$.

$Case\ 3.1.2.2.3.1$ $C(vw)=a,C(vz)=b\not=a,C(wx)=b,C(zx)=C(xy)=a,C(zy)=b,C(uw)=c,c\not=a,c\not=b$, there exists 
a $(wy,b)$-path, there exists a $(wy,a)$-path which does not go through the vertex $x$
and there exists a $(zu,a)$-path which does not go through the vertex $x$.
In this case we recolor $C(uw)=a,C(uv)=c,C(vw)=b,C(vz)=a,C(zx)=b,C(wx)=c$ as in Figure~\ref{fig:two_cubic_3_1_1_2_3_1}.
We do not introduce any monochromatic cycle because colors $a,b,c$ are pairwise different
and all the $a$-edges of the configuration are on an $a$-path which ends at $v$ and $x$.

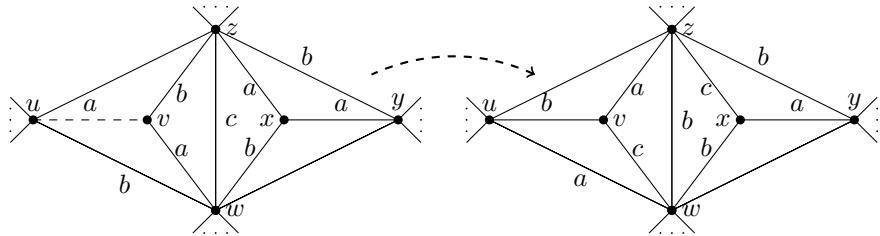
\begin{figure}[hptb]
	\begin{center}
  {\small{
  \begin{tikzpicture}[scale=0.6]

        \foreach \x/\y/\z in {0/0, 0/4, -1.5/2, -4/2}
    {
      \fill (\x,\y) circle (0.1);
      \fill (-\x,\y) circle (0.1);
    }

    \foreach \x in {-1, 1}
    {
      \draw (\x*4,2) -- (0,0) -- (0,4) -- (\x * 1.5, 2) -- (0,0) -- (\x * 4, 2) -- (0,4);
    }
    \draw (4,2) -- (1.5,2);
    \draw[dashed] (-4,2) -- (-1.5,2);

    \begin{scope}[shift={(0.0, 2.0)}]

      \foreach \x in {-1,1} 
      {
        \draw (0, \x * 2) -- (-0.5,\x * 2.5);
        \draw (0,\x * 2) -- (0.5, \x * 2.5);
        \fill (-0.25, \x * 2.5) circle (0.02);
        \fill (0,\x * 2.5) circle (0.02);
        \fill (0.25, \x * 2.5) circle (0.02);
      }

      \foreach \x in {-1,1} 
      {
        \draw (\x * 4, 0) -- (\x * 4.5, 0.5);
        \draw (\x * 4, 0) -- (\x * 4.5, -0.5);

        \fill (\x * 4.5, 0.25) circle (0.02);
        \fill (\x * 4.5, 0) circle (0.02);
        \fill (\x * 4.5, -0.25) circle (0.02);
      }
    \end{scope}

    \draw[above] (-4,2) node {$u$};
    \draw[above] (4,2) node {$y$};

    \draw[right] (0,4) node {$z$};
    \draw[right] (0,0) node {$w$};

    \draw[right] (-1.5,2) node {$v$};
    \draw[left] (1.5,2) node {$x$};

    \draw[above] (-2.75,2) node {$a$};
    \draw[above] (-0.75,1) node {$a$};
    \draw[below] (-0.75,3) node {$b$};
    \draw[above] (0.75,1) node {$b$};

    \draw[above] (2.75,2) node {$a$};
    \draw[below] (0.75,3) node {$a$};
    \draw[above] (2,3) node {$b$};
    \draw[below] (-2,1) node {$c$};

    \draw[dashed,thick,<-] (7,3) arc (65:115:3*1.414);

    \begin{scope}[shift={(10.0, 0.0)}]
          \foreach \x/\y/\z in {0/0, 0/4, -1.5/2, -4/2}
    {
      \fill (\x,\y) circle (0.1);
      \fill (-\x,\y) circle (0.1);
    }

    \foreach \x in {-1, 1}
    {
      \draw (\x*4,2) -- (0,0) -- (0,4) -- (\x * 1.5, 2) -- (0,0) -- (\x * 4, 2) -- (0,4);
    }
    \draw (4,2) -- (1.5,2);
    \draw (-4,2) -- (-1.5,2);

    \begin{scope}[shift={(0.0, 2.0)}]

      \foreach \x in {-1,1} 
      {
        \draw (0, \x * 2) -- (-0.5,\x * 2.5);
        \draw (0,\x * 2) -- (0.5, \x * 2.5);
        \fill (-0.25, \x * 2.5) circle (0.02);
        \fill (0,\x * 2.5) circle (0.02);
        \fill (0.25, \x * 2.5) circle (0.02);
      }

      \foreach \x in {-1,1} 
      {
        \draw (\x * 4, 0) -- (\x * 4.5, 0.5);
        \draw (\x * 4, 0) -- (\x * 4.5, -0.5);

        \fill (\x * 4.5, 0.25) circle (0.02);
        \fill (\x * 4.5, 0) circle (0.02);
        \fill (\x * 4.5, -0.25) circle (0.02);
      }
    \end{scope}

    \draw[above] (-4,2) node {$u$};
    \draw[above] (4,2) node {$y$};

    \draw[right] (0,4) node {$z$};
    \draw[right] (0,0) node {$w$};

    \draw[right] (-1.5,2) node {$v$};
    \draw[left] (1.5,2) node {$x$};

      \draw[above] (-2.75,2) node {$c$};
      \draw[above] (-0.75,1) node {$b$};
      \draw[below] (-0.75,3) node {$a$};
      \draw[above] (0.75,1) node {$c$};

      \draw[above] (2.75,2) node {$a$};
      \draw[below] (0.75,3) node {$b$};
      \draw[above] (2,3) node {$b$};
      \draw[below] (-2,1) node {$a$};
    \end{scope}

  \end{tikzpicture}
  }}
	\end{center}
	\caption{Case $3.1.2.2.3.1$ in the proof of Lemma~\ref{lem:two_cubic}}
	\label{fig:two_cubic_3_1_1_2_3_1}
\end{figure}

$Case\ 3.1.2.2.3.2$ $C(vw)=a,C(vz)=b\not=a,C(wx)=b,C(zx)=C(xy)=a,C(zy)=b,C(uw)=b$, there exists 
a $(uy,b)$-path, there exists a $(wy,a)$-path which does not go through the vertex $x$
and there exists a $(zu,a)$-path which does not go through the vertex $x$.
Let $c=C(wz)$.
Since the vertex $w$ has already two incident $b$-edges we have $c\not=b$.
Because of the $(wy,a)$-path we have $c\not=a$ thus we can recolor
$C(uw)=a,C(uv)=b,C(vw)=c,C(vz)=a,C(zx)=c,C(zw)=b$ as in Figure~\ref{fig:two_cubic_3_1_1_2_3_2}.
We do not introduce any monochromatic cycle because colors $a,b,c$ are pairwise distinct
and all the $a$-edges of the configuration are on an $a$-path which ends at $v$ and $x$.

\begin{figure}[hptb]
	\begin{center}
  {\small{
  \begin{tikzpicture}[scale=0.6]

        \foreach \x/\y/\z in {0/0, 0/4, -1.5/2, -4/2}
    {
      \fill (\x,\y) circle (0.1);
      \fill (-\x,\y) circle (0.1);
    }

    \foreach \x in {-1, 1}
    {
      \draw (\x*4,2) -- (0,0) -- (0,4) -- (\x * 1.5, 2) -- (0,0) -- (\x * 4, 2) -- (0,4);
    }
    \draw (4,2) -- (1.5,2);
    \draw[dashed] (-4,2) -- (-1.5,2);

    \begin{scope}[shift={(0.0, 2.0)}]

      \foreach \x in {-1,1} 
      {
        \draw (0, \x * 2) -- (-0.5,\x * 2.5);
        \draw (0,\x * 2) -- (0.5, \x * 2.5);
        \fill (-0.25, \x * 2.5) circle (0.02);
        \fill (0,\x * 2.5) circle (0.02);
        \fill (0.25, \x * 2.5) circle (0.02);
      }

      \foreach \x in {-1,1} 
      {
        \draw (\x * 4, 0) -- (\x * 4.5, 0.5);
        \draw (\x * 4, 0) -- (\x * 4.5, -0.5);

        \fill (\x * 4.5, 0.25) circle (0.02);
        \fill (\x * 4.5, 0) circle (0.02);
        \fill (\x * 4.5, -0.25) circle (0.02);
      }
    \end{scope}

    \draw[above] (-4,2) node {$u$};
    \draw[above] (4,2) node {$y$};

    \draw[right] (0,4) node {$z$};
    \draw[right] (0,0) node {$w$};

    \draw[right] (-1.5,2) node {$v$};
    \draw[left] (1.5,2) node {$x$};

    \draw[above] (-2.75,2) node {$a$};
    \draw[above] (-0.75,1) node {$a$};
    \draw[below] (-0.75,3) node {$b$};
    \draw[above] (0.75,1) node {$b$};

    \draw[above] (2.75,2) node {$a$};
    \draw[below] (0.75,3) node {$a$};
    \draw[above] (2,3) node {$b$};
    \draw[below] (-2,1) node {$b$};
    \draw[right] (0,2) node {$c$};

    \draw[dashed,thick,<-] (7,3) arc (65:115:3*1.414);

    \begin{scope}[shift={(10.0, 0.0)}]
          \foreach \x/\y/\z in {0/0, 0/4, -1.5/2, -4/2}
    {
      \fill (\x,\y) circle (0.1);
      \fill (-\x,\y) circle (0.1);
    }

    \foreach \x in {-1, 1}
    {
      \draw (\x*4,2) -- (0,0) -- (0,4) -- (\x * 1.5, 2) -- (0,0) -- (\x * 4, 2) -- (0,4);
    }
    \draw (4,2) -- (1.5,2);
    \draw (-4,2) -- (-1.5,2);

    \begin{scope}[shift={(0.0, 2.0)}]

      \foreach \x in {-1,1} 
      {
        \draw (0, \x * 2) -- (-0.5,\x * 2.5);
        \draw (0,\x * 2) -- (0.5, \x * 2.5);
        \fill (-0.25, \x * 2.5) circle (0.02);
        \fill (0,\x * 2.5) circle (0.02);
        \fill (0.25, \x * 2.5) circle (0.02);
      }

      \foreach \x in {-1,1} 
      {
        \draw (\x * 4, 0) -- (\x * 4.5, 0.5);
        \draw (\x * 4, 0) -- (\x * 4.5, -0.5);

        \fill (\x * 4.5, 0.25) circle (0.02);
        \fill (\x * 4.5, 0) circle (0.02);
        \fill (\x * 4.5, -0.25) circle (0.02);
      }
    \end{scope}

    \draw[above] (-4,2) node {$u$};
    \draw[above] (4,2) node {$y$};

    \draw[right] (0,4) node {$z$};
    \draw[right] (0,0) node {$w$};

    \draw[right] (-1.5,2) node {$v$};
    \draw[left] (1.5,2) node {$x$};

      \draw[above] (-2.75,2) node {$b$};
      \draw[above] (-0.75,1) node {$c$};
      \draw[below] (-0.75,3) node {$a$};
      \draw[above] (0.75,1) node {$b$};

      \draw[above] (2.75,2) node {$a$};
      \draw[below] (0.75,3) node {$c$};
      \draw[above] (2,3) node {$b$};
      \draw[below] (-2,1) node {$a$};
      \draw[right] (0,2) node {$b$};

    \end{scope}

  \end{tikzpicture}
  }}
	\end{center}
	\caption{Case $3.1.2.2.3.2$ in the proof of Lemma~\ref{lem:two_cubic}}
	\label{fig:two_cubic_3_1_1_2_3_2}
\end{figure}

$Case\ 3.2$ $C(vw)=a,C(vz)=b\not=a,C(wx)=a$ and there exists a $(xu,a)$-path, see Figure~\ref{fig:two_cubic_3_2}.
Because of the $(xu,a)$-path the vertex $x$ has exactly two incident edges colored $a$ thus we consider
which is the second one.

\begin{figure}[hptb]
	\begin{center}
  {\small{
  \begin{tikzpicture}[scale=0.6]

        \foreach \x/\y/\z in {0/0, 0/4, -1.5/2, -4/2}
    {
      \fill (\x,\y) circle (0.1);
      \fill (-\x,\y) circle (0.1);
    }

    \foreach \x in {-1, 1}
    {
      \draw (\x*4,2) -- (0,0) -- (0,4) -- (\x * 1.5, 2) -- (0,0) -- (\x * 4, 2) -- (0,4);
    }
    \draw (4,2) -- (1.5,2);
    \draw[dashed] (-4,2) -- (-1.5,2);

    \begin{scope}[shift={(0.0, 2.0)}]

      \foreach \x in {-1,1} 
      {
        \draw (0, \x * 2) -- (-0.5,\x * 2.5);
        \draw (0,\x * 2) -- (0.5, \x * 2.5);
        \fill (-0.25, \x * 2.5) circle (0.02);
        \fill (0,\x * 2.5) circle (0.02);
        \fill (0.25, \x * 2.5) circle (0.02);
      }

      \foreach \x in {-1,1} 
      {
        \draw (\x * 4, 0) -- (\x * 4.5, 0.5);
        \draw (\x * 4, 0) -- (\x * 4.5, -0.5);

        \fill (\x * 4.5, 0.25) circle (0.02);
        \fill (\x * 4.5, 0) circle (0.02);
        \fill (\x * 4.5, -0.25) circle (0.02);
      }
    \end{scope}

    \draw[above] (-4,2) node {$u$};
    \draw[above] (4,2) node {$y$};

    \draw[right] (0,4) node {$z$};
    \draw[right] (0,0) node {$w$};

    \draw[right] (-1.5,2) node {$v$};
    \draw[left] (1.5,2) node {$x$};

    \draw[above] (-2.75,2) node {$a$};
    \draw[above] (-0.75,1) node {$a$};
    \draw[below] (-0.75,3) node {$b$};
    \draw[above] (0.75,1) node {$a$};

  \end{tikzpicture}
  }}
	\end{center}
	\caption{Case $3.2$ in the proof of Lemma~\ref{lem:two_cubic}}
	\label{fig:two_cubic_3_2}
\end{figure}

$Case\ 3.2.1$ $C(vw)=a,C(vz)=b\not=a,C(wx)=a,C(xy)=a$ and there exists a $(uy,a)$-path. 
Let $c=C(xz)$. Obviously $c\not=a$ since otherwise the vertex $x$ would have three incident
$a$-edges. Now we see what happens if $c\ne b$ and what if $c=b$.

$Case\ 3.2.1.1$ $C(vw)=a,C(vz)=b,C(wx)=a,C(xy)=a,C(xz)=c,a\not=b\not=c\not=a$ and there exists a $(uy,a)$-path. 
Note that we can assume that $C(wy)\ne c$, for otherwise we just swap the colors of $zv$ and $zx$.
Then we recolor edges $wx$ and $xy$ to $C(wy)$ and $wy$ to $a$ and we obtain the situation where $C(wx)\not=a$ which was already considered in Case $3.1$.

$Case\ 3.2.1.2$ $C(vw)=a,C(vz)=b,C(wx)=a,C(xy)=a,C(xz)=b,a\not=b$ and there exists a $(uy,a)$-path. 
In this case we use the color $c=C(uw)\not=a$ to recolor $C(xw)=C(uv)=c$ and $C(uw)=a$
as in Figure~\ref{fig:two_cubic_3_2_1_2}. Since $a\not=b$ and $a\not=c$ we do not introduce any $a$-cycle.
If $b\not=c$ we are done but if $b=c$ we can have a $b$-cycle if there exists a $(uw,b)$-path but
it can not be true because it would imply that there is a $b$-cycle in $G'$.

\begin{figure}[hptb]
	\begin{center}
  {\small{
  \begin{tikzpicture}[scale=0.6]

        \foreach \x/\y/\z in {0/0, 0/4, -1.5/2, -4/2}
    {
      \fill (\x,\y) circle (0.1);
      \fill (-\x,\y) circle (0.1);
    }

    \foreach \x in {-1, 1}
    {
      \draw (\x*4,2) -- (0,0) -- (0,4) -- (\x * 1.5, 2) -- (0,0) -- (\x * 4, 2) -- (0,4);
    }
    \draw (4,2) -- (1.5,2);
    \draw[dashed] (-4,2) -- (-1.5,2);

    \begin{scope}[shift={(0.0, 2.0)}]

      \foreach \x in {-1,1} 
      {
        \draw (0, \x * 2) -- (-0.5,\x * 2.5);
        \draw (0,\x * 2) -- (0.5, \x * 2.5);
        \fill (-0.25, \x * 2.5) circle (0.02);
        \fill (0,\x * 2.5) circle (0.02);
        \fill (0.25, \x * 2.5) circle (0.02);
      }

      \foreach \x in {-1,1} 
      {
        \draw (\x * 4, 0) -- (\x * 4.5, 0.5);
        \draw (\x * 4, 0) -- (\x * 4.5, -0.5);

        \fill (\x * 4.5, 0.25) circle (0.02);
        \fill (\x * 4.5, 0) circle (0.02);
        \fill (\x * 4.5, -0.25) circle (0.02);
      }
    \end{scope}

    \draw[above] (-4,2) node {$u$};
    \draw[above] (4,2) node {$y$};

    \draw[right] (0,4) node {$z$};
    \draw[right] (0,0) node {$w$};

    \draw[right] (-1.5,2) node {$v$};
    \draw[left] (1.5,2) node {$x$};

    \draw[above] (-2.75,2) node {$a$};
    \draw[above] (-0.75,1) node {$a$};
    \draw[below] (-0.75,3) node {$b$};
    \draw[below] (0.75,3) node {$b$};
    \draw[above] (0.75,1) node {$a$};
    \draw[above] (2.75,2) node {$a$};
    \draw[below] (-2,1) node {$c$};

    \draw[dashed,thick,<-] (7,3) arc (65:115:3*1.414);

    \begin{scope}[shift={(10.0, 0.0)}]

          \foreach \x/\y/\z in {0/0, 0/4, -1.5/2, -4/2}
    {
      \fill (\x,\y) circle (0.1);
      \fill (-\x,\y) circle (0.1);
    }

    \foreach \x in {-1, 1}
    {
      \draw (\x*4,2) -- (0,0) -- (0,4) -- (\x * 1.5, 2) -- (0,0) -- (\x * 4, 2) -- (0,4);
    }
    \draw (4,2) -- (1.5,2);
    \draw (-4,2) -- (-1.5,2);

    \begin{scope}[shift={(0.0, 2.0)}]

      \foreach \x in {-1,1} 
      {
        \draw (0, \x * 2) -- (-0.5,\x * 2.5);
        \draw (0,\x * 2) -- (0.5, \x * 2.5);
        \fill (-0.25, \x * 2.5) circle (0.02);
        \fill (0,\x * 2.5) circle (0.02);
        \fill (0.25, \x * 2.5) circle (0.02);
      }

      \foreach \x in {-1,1} 
      {
        \draw (\x * 4, 0) -- (\x * 4.5, 0.5);
        \draw (\x * 4, 0) -- (\x * 4.5, -0.5);

        \fill (\x * 4.5, 0.25) circle (0.02);
        \fill (\x * 4.5, 0) circle (0.02);
        \fill (\x * 4.5, -0.25) circle (0.02);
      }
    \end{scope}

    \draw[above] (-4,2) node {$u$};
    \draw[above] (4,2) node {$y$};

    \draw[right] (0,4) node {$z$};
    \draw[right] (0,0) node {$w$};

    \draw[right] (-1.5,2) node {$v$};
    \draw[left] (1.5,2) node {$x$};

      \draw[above] (-2.75,2) node {$c$};
      \draw[above] (-0.75,1) node {$a$};
      \draw[below] (-0.75,3) node {$b$};
      \draw[below] (0.75,3) node {$b$};
      \draw[above] (0.75,1) node {$c$};
      \draw[above] (2.75,2) node {$a$};
      \draw[below] (-2,1) node {$a$};

    \end{scope}

  \end{tikzpicture}
  }}
	\end{center}
	\caption{Case $3.2.1.2$ in the proof of Lemma~\ref{lem:two_cubic}}
	\label{fig:two_cubic_3_2_1_2}
\end{figure}

$Case\ 3.2.2$ $C(vw)=a,C(vz)=b\not=a,C(wx)=a,C(xz)=a$ and there exists a $(uz,a)$-path. 
Let $c=C(wz)$. We may assume that $C(xy)=c$ since otherwise we can recolor $C(wx)=C(xz)=c$
and $C(wz)=a$ ending in a situation where $C(wx)\not=a$ which is considered in Case $3.1$.
We swap colors of edges $zv$ and $zx$ (as in Figure~\ref{fig:two_cubic_3_2_2})
and check whether we introduce a monochromatic cycle (without taking the edge $uv$ into consideration).

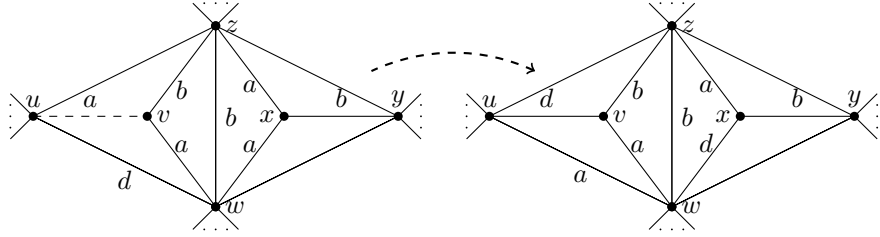
\begin{figure}[hptb]
	\begin{center}
  {\small{
  \begin{tikzpicture}[scale=0.6]

        \foreach \x/\y/\z in {0/0, 0/4, -1.5/2, -4/2}
    {
      \fill (\x,\y) circle (0.1);
      \fill (-\x,\y) circle (0.1);
    }

    \foreach \x in {-1, 1}
    {
      \draw (\x*4,2) -- (0,0) -- (0,4) -- (\x * 1.5, 2) -- (0,0) -- (\x * 4, 2) -- (0,4);
    }
    \draw (4,2) -- (1.5,2);
    \draw[dashed] (-4,2) -- (-1.5,2);

    \begin{scope}[shift={(0.0, 2.0)}]

      \foreach \x in {-1,1} 
      {
        \draw (0, \x * 2) -- (-0.5,\x * 2.5);
        \draw (0,\x * 2) -- (0.5, \x * 2.5);
        \fill (-0.25, \x * 2.5) circle (0.02);
        \fill (0,\x * 2.5) circle (0.02);
        \fill (0.25, \x * 2.5) circle (0.02);
      }

      \foreach \x in {-1,1} 
      {
        \draw (\x * 4, 0) -- (\x * 4.5, 0.5);
        \draw (\x * 4, 0) -- (\x * 4.5, -0.5);

        \fill (\x * 4.5, 0.25) circle (0.02);
        \fill (\x * 4.5, 0) circle (0.02);
        \fill (\x * 4.5, -0.25) circle (0.02);
      }
    \end{scope}

    \draw[above] (-4,2) node {$u$};
    \draw[above] (4,2) node {$y$};

    \draw[right] (0,4) node {$z$};
    \draw[right] (0,0) node {$w$};

    \draw[right] (-1.5,2) node {$v$};
    \draw[left] (1.5,2) node {$x$};

    \draw[above] (-2.75,2) node {$a$};
    \draw[above] (-0.75,1) node {$a$};
    \draw[below] (-0.75,3) node {$b$};

    \draw[right] (0,2) node {$c$};

    \draw[below] (0.75,3) node {$a$};
    \draw[above] (0.75,1) node {$a$};
    \draw[above] (2.75,2) node {$c$};

    \draw[dashed,thick,<-] (7,3) arc (65:115:3*1.414);

    \begin{scope}[shift={(10.0, 0.0)}]

          \foreach \x/\y/\z in {0/0, 0/4, -1.5/2, -4/2}
    {
      \fill (\x,\y) circle (0.1);
      \fill (-\x,\y) circle (0.1);
    }

    \foreach \x in {-1, 1}
    {
      \draw (\x*4,2) -- (0,0) -- (0,4) -- (\x * 1.5, 2) -- (0,0) -- (\x * 4, 2) -- (0,4);
    }
    \draw (4,2) -- (1.5,2);
    \draw[dashed] (-4,2) -- (-1.5,2);

    \begin{scope}[shift={(0.0, 2.0)}]

      \foreach \x in {-1,1} 
      {
        \draw (0, \x * 2) -- (-0.5,\x * 2.5);
        \draw (0,\x * 2) -- (0.5, \x * 2.5);
        \fill (-0.25, \x * 2.5) circle (0.02);
        \fill (0,\x * 2.5) circle (0.02);
        \fill (0.25, \x * 2.5) circle (0.02);
      }

      \foreach \x in {-1,1} 
      {
        \draw (\x * 4, 0) -- (\x * 4.5, 0.5);
        \draw (\x * 4, 0) -- (\x * 4.5, -0.5);

        \fill (\x * 4.5, 0.25) circle (0.02);
        \fill (\x * 4.5, 0) circle (0.02);
        \fill (\x * 4.5, -0.25) circle (0.02);
      }
    \end{scope}

    \draw[above] (-4,2) node {$u$};
    \draw[above] (4,2) node {$y$};

    \draw[right] (0,4) node {$z$};
    \draw[right] (0,0) node {$w$};

    \draw[right] (-1.5,2) node {$v$};
    \draw[left] (1.5,2) node {$x$};

      \draw[above] (-2.75,2) node {$a$};
      \draw[above] (-0.75,1) node {$a$};
      \draw[below] (-0.75,3) node {$a$};

      \draw[right] (0,2) node {$c$};

      \draw[below] (0.75,3) node {$b$};
      \draw[above] (0.75,1) node {$a$};
      \draw[above] (2.75,2) node {$c$};

    \end{scope}

  \end{tikzpicture}
  }}
	\end{center}
	\caption{Case $3.2.2$ in the proof of Lemma~\ref{lem:two_cubic}}
	\label{fig:two_cubic_3_2_2}
\end{figure}

$Case\ 3.2.2.1$ We do not introduce a monochromatic cycle thus we have a $k$-linear coloring
satisfying $C(wv)=C(vz)=a$ which is considered in Case $2$.

$Case\ 3.2.2.2$ We introduce a monochromatic cycle. The only possibility which leads to a monochromatic
cycle is when $b=c$ and there exists a $(wy,b)$-path. Now let us take into consideration edge $uw$ colored $d=C(uw)$.
We try to recolor $C(uv)=C(wx)=d$ and $C(uw)=a$ as in Figure~\ref{fig:two_cubic_3_2_2_2}.
If there is no monochromatic cycle even when taking the edge $uv$ into consideration we are done.
Otherwise we can assume that $b=d$ and we consider subcases regarding the color of the edge $uz$.

\begin{figure}[hptb]
	\begin{center}
  {\small{
  \begin{tikzpicture}[scale=0.6]

        \foreach \x/\y/\z in {0/0, 0/4, -1.5/2, -4/2}
    {
      \fill (\x,\y) circle (0.1);
      \fill (-\x,\y) circle (0.1);
    }

    \foreach \x in {-1, 1}
    {
      \draw (\x*4,2) -- (0,0) -- (0,4) -- (\x * 1.5, 2) -- (0,0) -- (\x * 4, 2) -- (0,4);
    }
    \draw (4,2) -- (1.5,2);
    \draw[dashed] (-4,2) -- (-1.5,2);

    \begin{scope}[shift={(0.0, 2.0)}]

      \foreach \x in {-1,1} 
      {
        \draw (0, \x * 2) -- (-0.5,\x * 2.5);
        \draw (0,\x * 2) -- (0.5, \x * 2.5);
        \fill (-0.25, \x * 2.5) circle (0.02);
        \fill (0,\x * 2.5) circle (0.02);
        \fill (0.25, \x * 2.5) circle (0.02);
      }

      \foreach \x in {-1,1} 
      {
        \draw (\x * 4, 0) -- (\x * 4.5, 0.5);
        \draw (\x * 4, 0) -- (\x * 4.5, -0.5);

        \fill (\x * 4.5, 0.25) circle (0.02);
        \fill (\x * 4.5, 0) circle (0.02);
        \fill (\x * 4.5, -0.25) circle (0.02);
      }
    \end{scope}

    \draw[above] (-4,2) node {$u$};
    \draw[above] (4,2) node {$y$};

    \draw[right] (0,4) node {$z$};
    \draw[right] (0,0) node {$w$};

    \draw[right] (-1.5,2) node {$v$};
    \draw[left] (1.5,2) node {$x$};

    \draw[above] (-2.75,2) node {$a$};
    \draw[above] (-0.75,1) node {$a$};
    \draw[below] (-0.75,3) node {$b$};

    \draw[right] (0,2) node {$b$};

    \draw[below] (0.75,3) node {$a$};
    \draw[above] (0.75,1) node {$a$};
    \draw[above] (2.75,2) node {$b$};

    \draw[below] (-2,1) node {$d$};

    \draw[dashed,thick,<-] (7,3) arc (65:115:3*1.414);

    \begin{scope}[shift={(10.0, 0.0)}]

          \foreach \x/\y/\z in {0/0, 0/4, -1.5/2, -4/2}
    {
      \fill (\x,\y) circle (0.1);
      \fill (-\x,\y) circle (0.1);
    }

    \foreach \x in {-1, 1}
    {
      \draw (\x*4,2) -- (0,0) -- (0,4) -- (\x * 1.5, 2) -- (0,0) -- (\x * 4, 2) -- (0,4);
    }
    \draw (4,2) -- (1.5,2);
    \draw (-4,2) -- (-1.5,2);

    \begin{scope}[shift={(0.0, 2.0)}]

      \foreach \x in {-1,1} 
      {
        \draw (0, \x * 2) -- (-0.5,\x * 2.5);
        \draw (0,\x * 2) -- (0.5, \x * 2.5);
        \fill (-0.25, \x * 2.5) circle (0.02);
        \fill (0,\x * 2.5) circle (0.02);
        \fill (0.25, \x * 2.5) circle (0.02);
      }

      \foreach \x in {-1,1} 
      {
        \draw (\x * 4, 0) -- (\x * 4.5, 0.5);
        \draw (\x * 4, 0) -- (\x * 4.5, -0.5);

        \fill (\x * 4.5, 0.25) circle (0.02);
        \fill (\x * 4.5, 0) circle (0.02);
        \fill (\x * 4.5, -0.25) circle (0.02);
      }
    \end{scope}

    \draw[above] (-4,2) node {$u$};
    \draw[above] (4,2) node {$y$};

    \draw[right] (0,4) node {$z$};
    \draw[right] (0,0) node {$w$};

    \draw[right] (-1.5,2) node {$v$};
    \draw[left] (1.5,2) node {$x$};

      \draw[above] (-2.75,2) node {$d$};
      \draw[above] (-0.75,1) node {$a$};
      \draw[below] (-0.75,3) node {$b$};

      \draw[right] (0,2) node {$b$};

      \draw[below] (0.75,3) node {$a$};
      \draw[above] (0.75,1) node {$d$};
      \draw[above] (2.75,2) node {$b$};

      \draw[below] (-2,1) node {$a$};

    \end{scope}

  \end{tikzpicture}
  }}
	\end{center}
	\caption{Case $3.2.2.2$ in the proof of Lemma~\ref{lem:two_cubic}}
	\label{fig:two_cubic_3_2_2_2}
\end{figure}

$Case\ 3.2.2.2.1$ $C(vw)=C(wx)=C(xz)=a,C(xy)=C(vz)=C(zw)=C(wu)=b\not=a,C(uz)=a$. Let us
consider the edge $wy$ colored $c=C(wy)$. Since the vertex $w$ already has two incident
$a$-edges and two incident $b$-edges we have $c\not=a$ and $c\not=b$.
In this case we recolor $C(uz)=b,C(uw)=a,C(vw)=c,C(wz)=a,C(wx)=b,C(wy)=b,C(xy)=c$ and put $C(uv)=a$
as in Figure~\ref{fig:two_cubic_3_2_2_2_1} without introducing any monochromatic cycle.

\begin{figure}[hptb]
	\begin{center}
  {\small{
  \begin{tikzpicture}[scale=0.6]

        \foreach \x/\y/\z in {0/0, 0/4, -1.5/2, -4/2}
    {
      \fill (\x,\y) circle (0.1);
      \fill (-\x,\y) circle (0.1);
    }

    \foreach \x in {-1, 1}
    {
      \draw (\x*4,2) -- (0,0) -- (0,4) -- (\x * 1.5, 2) -- (0,0) -- (\x * 4, 2) -- (0,4);
    }
    \draw (4,2) -- (1.5,2);
    \draw[dashed] (-4,2) -- (-1.5,2);

    \begin{scope}[shift={(0.0, 2.0)}]

      \foreach \x in {-1,1} 
      {
        \draw (0, \x * 2) -- (-0.5,\x * 2.5);
        \draw (0,\x * 2) -- (0.5, \x * 2.5);
        \fill (-0.25, \x * 2.5) circle (0.02);
        \fill (0,\x * 2.5) circle (0.02);
        \fill (0.25, \x * 2.5) circle (0.02);
      }

      \foreach \x in {-1,1} 
      {
        \draw (\x * 4, 0) -- (\x * 4.5, 0.5);
        \draw (\x * 4, 0) -- (\x * 4.5, -0.5);

        \fill (\x * 4.5, 0.25) circle (0.02);
        \fill (\x * 4.5, 0) circle (0.02);
        \fill (\x * 4.5, -0.25) circle (0.02);
      }
    \end{scope}

    \draw[above] (-4,2) node {$u$};
    \draw[above] (4,2) node {$y$};

    \draw[right] (0,4) node {$z$};
    \draw[right] (0,0) node {$w$};

    \draw[right] (-1.5,2) node {$v$};
    \draw[left] (1.5,2) node {$x$};

    \draw[above] (-2.75,2) node {$a$};
    \draw[above] (-0.75,1) node {$a$};
    \draw[below] (-0.75,3) node {$b$};

    \draw[right] (0,2) node {$b$};

    \draw[below] (0.75,3) node {$a$};
    \draw[above] (0.75,1) node {$a$};
    \draw[above] (2.75,2) node {$b$};

    \draw[below] (-2,1) node {$b$};
    \draw[above] (-2,3) node {$a$};

    \draw[below] (2,1) node {$c$};

    \draw[dashed,thick,<-] (7,3) arc (65:115:3*1.414);

    \begin{scope}[shift={(10.0, 0.0)}]

          \foreach \x/\y/\z in {0/0, 0/4, -1.5/2, -4/2}
    {
      \fill (\x,\y) circle (0.1);
      \fill (-\x,\y) circle (0.1);
    }

    \foreach \x in {-1, 1}
    {
      \draw (\x*4,2) -- (0,0) -- (0,4) -- (\x * 1.5, 2) -- (0,0) -- (\x * 4, 2) -- (0,4);
    }
    \draw (4,2) -- (1.5,2);
    \draw (-4,2) -- (-1.5,2);

    \begin{scope}[shift={(0.0, 2.0)}]

      \foreach \x in {-1,1} 
      {
        \draw (0, \x * 2) -- (-0.5,\x * 2.5);
        \draw (0,\x * 2) -- (0.5, \x * 2.5);
        \fill (-0.25, \x * 2.5) circle (0.02);
        \fill (0,\x * 2.5) circle (0.02);
        \fill (0.25, \x * 2.5) circle (0.02);
      }

      \foreach \x in {-1,1} 
      {
        \draw (\x * 4, 0) -- (\x * 4.5, 0.5);
        \draw (\x * 4, 0) -- (\x * 4.5, -0.5);

        \fill (\x * 4.5, 0.25) circle (0.02);
        \fill (\x * 4.5, 0) circle (0.02);
        \fill (\x * 4.5, -0.25) circle (0.02);
      }
    \end{scope}

    \draw[above] (-4,2) node {$u$};
    \draw[above] (4,2) node {$y$};

    \draw[right] (0,4) node {$z$};
    \draw[right] (0,0) node {$w$};

    \draw[right] (-1.5,2) node {$v$};
    \draw[left] (1.5,2) node {$x$};

      \draw[above] (-2.75,2) node {$a$};
      \draw[above] (-0.75,1) node {$c$};
      \draw[below] (-0.75,3) node {$b$};

      \draw[right] (0,2) node {$a$};

      \draw[below] (0.75,3) node {$a$};
      \draw[above] (0.75,1) node {$b$};
      \draw[above] (2.75,2) node {$c$};

      \draw[below] (-2,1) node {$a$};
      \draw[above] (-2,3) node {$b$};

      \draw[below] (2,1) node {$b$};

    \end{scope}

  \end{tikzpicture}
  }}
	\end{center}
	\caption{Case $3.2.2.2.1$ in the proof of Lemma~\ref{lem:two_cubic}}
	\label{fig:two_cubic_3_2_2_2_1}
\end{figure}

$Case\ 3.2.2.2.2$ $C(vw)=C(wx)=C(xz)=a,C(xy)=C(vz)=C(zw)=C(wu)=b\not=a,C(uz)=c\not=a$.
Since the vertex $z$ already has two incident $b$ edges we have $c\not=b$.
In this case we recolor $C(uz)=b,C(zv)=a,C(zx)=c,C(uw)=a,C(vw)=b$ and put $C(uv)=c$ as in Figure~\ref{fig:two_cubic_3_2_2_2_2}.
We do not introduce any monochromatic cycle since 
each of the recolored edges is on a monochromatic path which ends at $v$ or $x$ (see Fig.~\ref{fig:two_cubic_3_2_2_2_2}).

\begin{figure}[hptb]
	\begin{center}
  {\small{
  \begin{tikzpicture}[scale=0.6]

        \foreach \x/\y/\z in {0/0, 0/4, -1.5/2, -4/2}
    {
      \fill (\x,\y) circle (0.1);
      \fill (-\x,\y) circle (0.1);
    }

    \foreach \x in {-1, 1}
    {
      \draw (\x*4,2) -- (0,0) -- (0,4) -- (\x * 1.5, 2) -- (0,0) -- (\x * 4, 2) -- (0,4);
    }
    \draw (4,2) -- (1.5,2);
    \draw[dashed] (-4,2) -- (-1.5,2);

    \begin{scope}[shift={(0.0, 2.0)}]

      \foreach \x in {-1,1} 
      {
        \draw (0, \x * 2) -- (-0.5,\x * 2.5);
        \draw (0,\x * 2) -- (0.5, \x * 2.5);
        \fill (-0.25, \x * 2.5) circle (0.02);
        \fill (0,\x * 2.5) circle (0.02);
        \fill (0.25, \x * 2.5) circle (0.02);
      }

      \foreach \x in {-1,1} 
      {
        \draw (\x * 4, 0) -- (\x * 4.5, 0.5);
        \draw (\x * 4, 0) -- (\x * 4.5, -0.5);

        \fill (\x * 4.5, 0.25) circle (0.02);
        \fill (\x * 4.5, 0) circle (0.02);
        \fill (\x * 4.5, -0.25) circle (0.02);
      }
    \end{scope}

    \draw[above] (-4,2) node {$u$};
    \draw[above] (4,2) node {$y$};

    \draw[right] (0,4) node {$z$};
    \draw[right] (0,0) node {$w$};

    \draw[right] (-1.5,2) node {$v$};
    \draw[left] (1.5,2) node {$x$};

    \draw[above] (-2.75,2) node {$a$};
    \draw[above] (-0.75,1) node {$a$};
    \draw[below] (-0.75,3) node {$b$};

    \draw[right] (0,2) node {$b$};

    \draw[below] (0.75,3) node {$a$};
    \draw[above] (0.75,1) node {$a$};
    \draw[above] (2.75,2) node {$b$};

    \draw[below] (-2,1) node {$b$};
    \draw[above] (-2,3) node {$c$};

    \draw[dashed,thick,<-] (7,3) arc (65:115:3*1.414);

    \begin{scope}[shift={(10.0, 0.0)}]

          \foreach \x/\y/\z in {0/0, 0/4, -1.5/2, -4/2}
    {
      \fill (\x,\y) circle (0.1);
      \fill (-\x,\y) circle (0.1);
    }

    \foreach \x in {-1, 1}
    {
      \draw (\x*4,2) -- (0,0) -- (0,4) -- (\x * 1.5, 2) -- (0,0) -- (\x * 4, 2) -- (0,4);
    }
    \draw (4,2) -- (1.5,2);
    \draw (-4,2) -- (-1.5,2);

    \begin{scope}[shift={(0.0, 2.0)}]

      \foreach \x in {-1,1} 
      {
        \draw (0, \x * 2) -- (-0.5,\x * 2.5);
        \draw (0,\x * 2) -- (0.5, \x * 2.5);
        \fill (-0.25, \x * 2.5) circle (0.02);
        \fill (0,\x * 2.5) circle (0.02);
        \fill (0.25, \x * 2.5) circle (0.02);
      }

      \foreach \x in {-1,1} 
      {
        \draw (\x * 4, 0) -- (\x * 4.5, 0.5);
        \draw (\x * 4, 0) -- (\x * 4.5, -0.5);

        \fill (\x * 4.5, 0.25) circle (0.02);
        \fill (\x * 4.5, 0) circle (0.02);
        \fill (\x * 4.5, -0.25) circle (0.02);
      }
    \end{scope}

    \draw[above] (-4,2) node {$u$};
    \draw[above] (4,2) node {$y$};

    \draw[right] (0,4) node {$z$};
    \draw[right] (0,0) node {$w$};

    \draw[right] (-1.5,2) node {$v$};
    \draw[left] (1.5,2) node {$x$};

      \draw[above] (-2.75,2) node {$c$};
      \draw[above] (-0.75,1) node {$b$};
      \draw[below] (-0.75,3) node {$a$};

      \draw[right] (0,2) node {$b$};

      \draw[below] (0.75,3) node {$c$};
      \draw[above] (0.75,1) node {$a$};
      \draw[above] (2.75,2) node {$b$};

      \draw[below] (-2,1) node {$a$};
      \draw[above] (-2,3) node {$b$};

    \end{scope}

  \end{tikzpicture}
  }}
	\end{center}
	\caption{Case $3.2.2.2.2$ in the proof of Lemma~\ref{lem:two_cubic}}
	\label{fig:two_cubic_3_2_2_2_2}
\end{figure}

\end{proof}

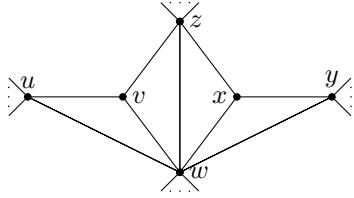
\begin{figure}[ht]
\centering
  {\small{  \begin{tikzpicture}[scale=0.5]
    \foreach \x/\y/\z in {0/0, 0/4, -1.5/2, -4/2}
    {
      \fill (\x,\y) circle (0.1);
      \fill (-\x,\y) circle (0.1);
    }

    \foreach \x in {-1, 1}
    {
      \draw (\x*4,2) -- (0,0) -- (0,4) -- (\x * 1.5, 2) -- (0,0) -- (\x * 4, 2);
    }
    \draw (4,2) -- (1.5,2);
    \draw (-4,2) -- (-1.5,2);

    \begin{scope}[shift={(0.0, 2.0)}]

      \foreach \x in {-1,1} 
      {
        \draw (0, \x * 2) -- (-0.5,\x * 2.5);
        \draw (0,\x * 2) -- (0.5, \x * 2.5);
        \fill (-0.25, \x * 2.5) circle (0.02);
        \fill (0,\x * 2.5) circle (0.02);
        \fill (0.25, \x * 2.5) circle (0.02);
      }

      \foreach \x in {-1,1} 
      {
        \draw (\x * 4, 0) -- (\x * 4.5, 0.5);
        \draw (\x * 4, 0) -- (\x * 4.5, -0.5);

        \fill (\x * 4.5, 0.25) circle (0.02);
        \fill (\x * 4.5, 0) circle (0.02);
        \fill (\x * 4.5, -0.25) circle (0.02);
      }
    \end{scope}

    \draw[above] (-4,2) node {$u$};
    \draw[above] (4,2) node {$y$};

    \draw[right] (0,4) node {$z$};
    \draw[right] (0,0) node {$w$};

    \draw[right] (-1.5,2) node {$v$};
    \draw[left] (1.5,2) node {$x$};

  \end{tikzpicture}
  }}
	\caption{The configuration from Corollary~\ref{cor:two_cubic}.}
	\label{fig:two_cubic_no_edges}
\end{figure}

\begin{corollary}
  \label{cor:two_cubic}
	$G$ does not contain the configuration in Fig.~\ref{fig:two_cubic_no_edges}.
\end{corollary}

\begin{proof}
From Lemma~\ref{lem:two-pairs} both $uz\in E(G)$ and $yz\in E(G)$. Hence we use Lemma~\ref{lem:two_cubic}. 
\end{proof}

\subsection{Proof of Proposition~\ref{prop:main}}

Now we use the following lemma due to Cole, Kowalik and \v{S}krekovski.

\begin{lemma}[Proposition 1.3 in~\cite{Kotzig}]
\label{lem:11-neighbors}
Let $G$ be a simple planar graph with minimum degree $\delta\geq 2$ such that each
$d$-vertex, $d\geq 12$, has at most $d-11$ neighbors of degree $2$. Then $G$ contains an edge of weight at most 13.
\end{lemma}

Using the above lemma and lemmas~\ref{lem_edge} and~\ref{lem:atmosttwo2vertices} only we can prove the following special case of
Proposition~\ref{prop:main}, which already improves known bounds on linear arboricity of planar graphs.

\begin{proposition}
\label{prop:special}
Any simple planar graph of maximum degree $\Delta$ has a linear coloring in $\max\{\ceil{\tfrac{\Delta}{2}}, 6\}$ colors.
\end{proposition}

\begin{proof}
Assume the claim is false and let $G$ be a minimal counterexample (in terms of the number of edges).
Let $k=\max\{\ceil{\tfrac{\Delta}{2}}, 6\}$. In particular, $\Delta \le 2k$.
By Lemma~\ref{lem_edge}, $G$ has no vertices of degree $1$ and any $2$-vertex has two neighbors of degree $2k$.
Next, by Lemma~\ref{lem:atmosttwo2vertices} every $(2k)$-vertex has at most one $2$-neighbor. Since $2k\ge 12$ the
assumptions of Lemma~\ref{lem:11-neighbors} are satisfied, so $G$ contains an edge of weight at most $13 \le 2k + 1$, a 
contradiction with Lemma~\ref{lem_edge}.
\end{proof}

Now, we proceed to the proof of Proposition~\ref{prop:main}. 
By the above proposition, Proposition~\ref{prop:main} holds for $\Delta\ge 11$. 
Hence, in what follows we assume that $\Delta \le 10$.
We put $k=\max\{\ceil{\tfrac{\Delta}{2}}, 5\}$ and we assume that $G$ is a minimal counterexample (in terms of the number of edges).

We prove Proposition~\ref{prop:main}  using the discharging method. The procedure is the following. 
We assign a number (called {\em charge}) to every vertex and face of a plane embedding of $G$, such that the total sum of all charges
is negative. Next, we redistribute the charge among vertices and faces in such a way, that
using the structural properties of graph $G$ described in Section~\ref{sec:structure} we are able to show that
every vertex and every face has nonnegative charge at the end, hence the total charge of
$G$ is nonnegative. This will give a contradiction, so the minimal counterexample does not
exist.

\paragraph{Initial charge.} We set the following initial charge to all vertices and faces of $G$:
\begin{eqnarray*}
	{\rm ch}_0(v) &=& \deg(v) - 4 \,,\,\,\, v \in V(G)\,,\\
	{\rm ch}_0(f) &=& \ell(f) - 4 \,,\,\,\, f \in F(G)\,.
\end{eqnarray*}
From Euler's formula we infer the following total charge of $G$:
\begin{eqnarray*}
	\sum_{v \in V(G)}(\deg(v)-4) + \displaystyle \sum_{f \in F(G)}(\ell(f)-4) = \\
	2|E(G)|-4|V(G)| + 2|E(G)| - 4 |F(G)| = \\
	-4 (|V(G)| - |E(G)| + |F(G)|) = -8. 
\end{eqnarray*}

\paragraph{Discharging rules.} Now, we present the discharging rules, by which we redistribute the charge
of vertices and faces in $G$. 

\begin{itemize}
	\item[(R1)] Every $10$-vertex sends $1$ to adjacent $2$-vertex.
	
	\item[(R2)] Every $(\ge 9)$-vertex sends $\tfrac{1}{3}$ to every adjacent $3$-vertex.

	\item[(R3)] Every $(\ge 8)$-vertex sends $\tfrac{1}{2}$ to every incident $3$-face with a vertex of degree at most 4.
				
	\item[(R4)] Every $(\ge 7)$-vertex sends $\tfrac{2}{5}$ to every incident $3$-face with a $5$-vertex.
				
	\item[(R5)] Every $(\ge 6)$-vertex sends $\tfrac{1}{3}$ to every incident $3$-face which is incident
				to only $(\ge\!6)$-vertices.
				
	\item[(R6)] Every $5$-vertex sends $\tfrac{1}{5}$ to every incident $3$-face.
	
	\item[(R7)] Every $(\ge \! 5)$-face $f$ sends $\tfrac{1}{3}$ to every incident $10$-vertex which has a $2$-neighbor incident to $f$.
\end{itemize}

\paragraph{Final charge.} 
Note that the initial charge is negative only for $2$- and $3$-vertices, and for $3$-faces.
	We show that by applying the discharging rules, all vertices and faces of $G$ have nonnegative final charge.
	
	First, we consider the charge of the faces. Note that $4$-faces do not send any charge so their charge remains $0$.
	Now we consider a face $f$ of length $\ell(f)\ge 5$. By Lemmas~\ref{lem_edge} and~\ref{lem:atmosttwo2vertices}, $f$ is incident 
	to at most $\floor{\tfrac{\ell(f)}{3}}$ vertices of degree 2. Hence $f$ sends at most $\tfrac{1}{3} \cdot 2 \cdot \lfloor \tfrac{\ell(f)}{3} \rfloor$ units
	of charge by (R7), which is less than $\ell(f)-4$ for $\ell(f)\ge 5$, hence $f$ retains positive charge.
	
	It only remains to show that every $3$-face $f$ receives at least $1$ from its neighbors, since its initial charge is $-1$.
	We consider cases regarding the degree of vertices incident to $f$. If $f$ is incident with a $2$-, $3$-,
	or $4$-vertex, it follows by Lemma~\ref{lem_edge} that the other two vertices incident with $f$ are of degree at least $2k-2\ge 8$,
	hence each of them sends $\tfrac{1}{2}$ to $f$ by (R3). 
	Next, if $f$ is incident to a $5$-vertex $v$, the other two incident vertices of $f$ are of degree at least $2k-3\ge 7$ by Lemma~\ref{lem_edge}. 
	Hence, $f$ receives $\tfrac{1}{5}$ from $v$ by (R6) and $\tfrac{2}{5}$ from each of the other two incident vertices by (R4), that is $1$ in total.
	Finally, if $f$ is incident only to $\ge\!6$-vertices, each of them sends $\tfrac{1}{3}$ by (R5), hence
	$f$ receives $1$ in total again. It follows that the final charge of $3$-faces is $0$.

\medskip
Now, we consider the final charge of the vertices. 
For convenience, we introduce a notion of a side. Let $v$ be a vertex and let $vx_0,\ldots,vx_{\deg(v)-1}$ be the edges incident to $v$,
enumerated in the clockwise order around $v$ in the given plane embedding. For any $i=0,\ldots,\deg(v)-1$, the pair $s=(vx_i,vx_{i+1})$ will be called a {\em side of $v$} (where $x_{\deg(v)}=x_0$). If $x_i$ and $x_{i+1}$ are adjacent, we say that $s$ is {\em triangular}.
We also say that $s$ is incident to $vx_i$ and $vx_{i+1}$.
Note that $v$ can have less than $\deg(v)$ incident faces (when $v$ is a cutvertex), while it has always $\deg(v)$ distinct incident sides.
However, for each triangular face incident to $v$ there is a distinct triangular side of $f$.
Since $v$ does not send charge to non-triangular faces, when $v$ sends charge to a triangle we can say that it sends the charge to the corresponding side and the total charge sent to sides is equal to the total charge sent to faces.
In what follows, we use the following claim.

\medskip

\noindent {\bf Claim 1} 
{\em If a $d$-vertex $v$ has negative final charge and $v$ is not adjacent to a $2$-vertex then $v$ has at most $11-d$ triangular sides.}

\bigskip
\noindent {\em Proof (of the claim).}
Let $p$ be the number of non-triangular sides of $v$.
Note that $v$ sends charge only to incident triangles (at most $\tfrac{1}{2}$ per triangle) and
to adjacent $3$-vertices (at most $\tfrac{1}{3}$ per $3$-vertex). For the proof of this claim, 
we replace (R2) by an equivalent rule:
\begin{itemize}
	\item[(R2')] For each $3$-neighbor $w$ of a $(\ge 9)$-vertex $v$, vertex $v$ sends $1/6$ to each of the two sides
incident with edge $vw$ and each of these sides resend the $1/6$ to $w$. 
\end{itemize}
Then, $v$ sends at most $\tfrac{2}{3}$ to each incident triangular side (the corresponding $3$-face has only one $3$-vertex, for otherwise
there is an edge of weight $6 < 2k+2$, which contradicts Lemma~\ref{lem_edge})
and it sends at most $\tfrac{1}{3}$ to each incident non-triangular side. It follows that $v$ sends 
at most $\tfrac{1}{3}\cdot p + \tfrac{2}{3}\cdot (d-p)=\tfrac{2d-p}{3}$ in total.
Hence, the final charge at $v$ is negative when $\tfrac{2d-p}{3} > d-4$, which is 
equivalent to $p < 12 - d$ and Claim 1 follows since $t$ is a natural number. \qed

Now we consider several cases regarding the degree of vertex $v$.
	\begin{itemize}
		\item{} {\it $v$ is a $2$-vertex.} The initial charge of $v$ is $-2$. 
		                By Lemma~\ref{lem_edge}, both its neighbors are of degree at least $2k\ge 10$.
				Hence, by (R1), $v$ receives $1$ from each of the two neighbors, and since it does not send any
				charge, its final charge is $0$. 

                \item{} {\it $v$ is a $3$-vertex.} The initial charge of $v$ is $-1$. 
		                By Lemma~\ref{lem_edge}, all three of its neighbors are of degree at least $2k-1\ge 9$.
				By (R2), $v$ receives $\tfrac{1}{3}$ from each of the three neighbors, and since it does not send any
				charge, its final charge is $0$. 
				
		\item{} {\it $v$ is a $4$-vertex.} In this case, $v$ does not send nor receive any charge. Hence, its initial
				charge, which is $0$, is equal to its final charge.
		
		\item{} {\it $v$ is a $d$-vertex, $5\le d \le 8$.} Note that $v$ sends charge only to incident triangles.
		                By rules (R3)-(R6), the charge $v$ sends to each incident triangle is at most 
				$\tfrac{1}{5}$, $\tfrac{1}{3}$, $\tfrac{2}{5}$, $\tfrac{1}{2}$, for $d=5,6,7,8$ respectively.
				One can check that in each of the four cases this is not more than $\tfrac{d-4}{d}$, and since there are at 
				most $d$ incident triangles, $v$ sends at most $d-4$ charge in total and retains nonnegative charge.
				
		\item{} {\it $v$ is a $9$-vertex.} The initial charge at $v$ is $5$.
		If $v$ has at most one $3$-neighbor then $v$ sends at most $\frac{9}{2}$ to faces and $\frac{1}{3}$ to vertices so
		its final charge is positive.
		If $v$ has at least two $3$-neighbors, then by Lemma~\ref{lem_cubicsmall}, each of them is incident with two non-triangular
		sides of $v$. Hence $v$ has at least 3 non-triangular sides, which contradicts Claim 1.

                \item{} {\it $v$ is a $10$-vertex.} 
		Assume first that $v$ has no $2$-neighbors. By Claim 1, $v$ is incident to at most one non-triangular side.
		If $v$ is incident only to $3$-faces, by Lemma~\ref{cor:two_cubic} $v$ has at most three $3$-neighbors.
		Then, $v$ sends at most $3\cdot\tfrac{1}{3}=1$ to vertices and at most $10\cdot\tfrac{1}{2}=5$ to faces,
		hence at most $6$ in total.
		If $v$ is incident to one non-triangular side, by Lemma~\ref{cor:two_cubic} $v$ has at most four $3$-neighbors.
		Then, $v$ sends at most $4\cdot\tfrac{1}{3}=\tfrac{4}{3}$ to vertices and at most $9\cdot\tfrac{1}{2}=\tfrac{9}{2}$ to faces,
		that is less than $6$ in total. Hence, in both cases the final charge at $v$ is nonnegative.
		
		Finally, assume that $v$ has a $2$-neighbor. Let $w_0, w_1, \ldots, w_9$ denote the neighbors of $v$ in the clockwise order
		in the given plane embedding of $G$. Assume w.l.o.g.\ $\deg(w_1)=2$.
		By Lemma~\ref{lem:atmosttwo2vertices}, this is the only $2$-neighbor of $v$.
		By Lemma~\ref{lem:2vertex}, the neighbors of $w_1$ are adjacent, so assume w.l.o.g.\ $w_0$ is adjacent with $w_1$ (in the beginning, we can choose the plane embedding of $G$ so that each triangle with one 2-vertex and two 10-vertices is a face).
		Since $G$ is simple, the face incident with $vw_1$ and $vw_2$, say $f$, is of length at least 4.
		
		Let $n_3$ denote the number of $3$-vertices among vertices $w_3,\ldots,w_9$.
		By Corollary~\ref{cor:atmostonetrianglewith3vertex}, each of these 3-neighbors is incident to at least one non-triangular side.
		Since each side is incident to at most two 3-neighbors of $v$, there are at least $\ceil{\frac{n_3}{2}}$ non-triangular sides, not counting the side $(vw_1,vw_2)$.
		
		It follows that $v$ sends 1 unit to $w_1$, $\frac{1}{3}$ to $w_2$ if $\deg(w_2)=3$,
		$\frac{n_3}{3}$ to vertices $w_3,\ldots,w_9$ and $\frac{1}{2}\cdot(9-\ceil{\frac{n_3}{2}})$ to incident faces.
		Hence, $v$ sends at most $5\frac{1}{2}+\frac{1}{3}[\deg(w_2)=3]+\frac{n_3}{3}-\frac{1}{2}\ceil{\frac{n_3}{2}}$.
		However, when $\deg(w_2)=3$, then face $f$ is of length at least 5 by Lemma~\ref{lem:23triangles}, so $v$ receives additional 
		$\tfrac{1}{3}$ from $f$ by (R7). Hence, $v$ gets at least $6+\frac{1}{3}[\deg(w_2)=3]$ charge in total. 
		It follows that the final charge at $v$ is at least $\frac{1}{2}-\frac{n_3}{3}+\frac{1}{2}\ceil{\frac{n_3}{2}}$, which is nonnegative 
		since $n_3\le 7$.
	\end{itemize}
It follows that the total charge of $G$ is nonnegative, establishing a contradiction on existence of a minimal counterexample.

\section{Algorithm}

In this section we show that our proof can be turned to an efficient algorithm for finding a linear $\lceil\tfrac{\Delta}{2}\rceil$-colorings.
The forbidden subgraphs from Section~\ref{sec:structure} will be called {\em reducible configurations}.
It should be clear the proof of Proposition~\ref{prop:main} corresponds to the following algorithm: find any of our reducible configurations in linear time,
then obtain a smaller graph in constant time by removing/contracting an edge, color it recursively and finally extend the coloring in linear time. Since in each recursive call the number of edges decreases, the number of recursive calls is linear, which gives $O(n^2)$ overall time complexity.
However, with some effort it is possible to improve the running time.
Namely, we present an $O(n\log n)$-time algorithm. The algorithm works for any planar graph and returns a partition into $\max\{\ceil{\tfrac{\Delta}{2}}, 5\}$ linear forests, which is optimal for $\Delta \ge 9$.

Our approach is as follows. First we describe an algorithm that finds a partition into $\max\{\ceil{\tfrac{\Delta}{2}}, 6\}$ linear forests, which is optimal for $\Delta \ge 11$. This can be treated an implementation of Proposition~\ref{prop:special}. Recall that for proving this proposition we needed only a few reducible configurations: an edge of weight at most $2k+1$, a 2-vertex with its neighbors nonadjacent, and a $2k$-vertex with two 2-neighbors.
As we will see in the Subsection~\ref{sec:alg-unbounded} these configurations are simple enough to find them very fast (even in constant time) after a linear preprocessing. 
Once we have this algorithm, we use it whenever $\Delta \ge 11$. Otherwise $\Delta \le 10$, so $\Delta$ is bounded which makes finding any bounded-size configuration very easy. Then we use the algorithm sketched in Subsection~\ref{sec:alg-bounded}.

\subsection{An algorithm for $\Delta \ge 11$}
\label{sec:alg-unbounded}

The coloring algorithm we describe in this section is inspired by the linear-time algorithm for $\Delta$-edge-coloring planar graphs presented in~\cite{edge_col}. 
For an input graph $G$ of maximum degree $\Delta$ we define $k=\max\{\ceil{\tfrac{\Delta}{2}}, 6\}$. 
We will show a $O(n\log n)$-time algorithm which finds a $k$-linear coloring of $G$.

Note that $\Delta \le 2k$, so the graph has vertices of degree a most $2k$ and $k\ge 6$
(we will use these facts in our arguments). We use the following three types of reducible edges of weight at most $2k+1$, which will be called \emph{nice}:

\begin{itemize}
\item edges of weight at most $13$,
\item edges incident to a $1$-vertex, and
\item edges incident to a $2$-vertex and a vertex of degree at most $2k-1$.
\end{itemize}

Our algorithm uses two queues: $Q_e$ and $Q_2$. Queue $Q_e$ stores nice edges, while queue $Q_2$ stores 2-vertices 
such that their both neighbors are of degree $2k$. 
Also, any $(2k)$-vertex $x$ may store a triangle $xyz$, such that $\deg(y)=2$ and $\deg(z)=2k$.
During the execution of the algorithm the following invariant is satisfied.

\begin{invariant}
\label{inv}
Any nice edge is stored in $Q_e$. 
Moreover, for any $2$-vertex $x$ with two $(2k)$-neighbors $v$ and $w$, either $x$ is $Q_2$ or $G$ contains a triangle $vxw$ and this triangle 
is stored in both $v$ and $w$.
Each vertex stores at most one triangle.
\end{invariant}

It is easy to initialize the queues in linear time to make Invariant~\ref{inv} satisfied at the beginning.
Then we use a recursive procedure which can be sketched as follows.
By configuration A and B we mean the configurations from the cases A and B of the proof of Lemma~\ref{lem:atmosttwo2vertices} (see Fig.~\ref{fig:two2vertices}).

\begin{enumerate}[Step 1.]
\item
  (Base of the recursion.) If $G$ has no edges, return the empty coloring.
\item
  If $Q_e$ contains an edge $e$, obtain a coloring of $G-e$ recursively and color $e$ by a free color as described in Lemma~\ref{lem_edge}.
\item
  Remove a 2-vertex $x$ from $Q_2$. Denote the neighbors of $x$ by $v$ and $w$. \label{step:remove-Q2}
\item
  If $v$ or $w$ stores a triangle $vwy$, we have configuration A. 
  Remove an edge $e$ of $G$ as described in Lemma~\ref{lem:atmosttwo2vertices}, recurse on $G-e$ and extend the coloring of $G-e$ to a coloring of $G$ as described in Lemma~\ref{lem:atmosttwo2vertices}.
\item 
  If $vw\not\in E(G)$, proceed as in Lemma~\ref{lem:2vertex}: remove vertex $x$ and add edge $vw$, recurse, add vertex $x$ and edges $vx$, $wx$, color these edges as $vw$ and remove edge $vw$. \label{step:contract}
\item 
  Else ($vw\in E(G)$)
  \begin{enumerate}[(i)]
   \item If $v$ (resp. $w$) stores a triangle $vyu$, we have configuration B. 
  Remove an edge $e$ of $G$ as described in Lemma~\ref{lem:atmosttwo2vertices}, recurse on $G-e$ and extend the coloring of $G-e$ to a coloring of $G$ as described in Lemma~\ref{lem:atmosttwo2vertices}.
   \item Otherwise, store triangle $vxw$ in $v$ and $w$. 
  \end{enumerate}
\end{enumerate}

Now we describe how the queues $Q_e$ and $Q_2$ are updated during an execution of the algorithm, to keep Invariant~\ref{inv} satisfied. 
First notice that it is easy to store degrees of vertices and update them in overall $O(n)$ time.
Then, whenever an edge is removed, for each of its endpoints,
say $z$, we check whether $z$ is of degree at most $12$. If so, for all its $O(1)$ incident edges
we check whether they are nice and if that is the case we add them to $Q_e$ (unless $Q_e$ already contains this edge).
Also, when after removing an edge a degree of its endpoint $z$ drops to 2, we check whether both of its neighbors are of degree $2k$ and if so,
we add $z$ to $Q_2$. Hence, updating $Q_e$ and $Q_2$ takes $O(1)$ time after each edge deletion. Clearly, after the graph modification in Step~\ref{step:contract}, there is no need to update any queue.

Now we are going to show the correctness of our algorithm. 

\begin{proposition}
Let $k=\max\{\ceil{\tfrac{\Delta}{2}}, 6\}$.
The above algorithm correctly finds a $k$-linear coloring of any planar graph of maximum degree $\Delta$.
\end{proposition}

\begin{proof}
Clearly, it suffices to show that whenever the algorithm finds itself in Step~\ref{step:remove-Q2} the queue $Q_2$ is not empty.
Assume the contrary. Since $Q_e$ is empty, i.e.\ there are no nice edges, $G$ has no 1-vertices and each 2-vertex is adjacent to two $(2k)$-vertices.
Hence only $(2k)$-vertices have $2$-neighbors. Since $Q_2$ is empty, by Invariant~\ref{inv} each $(2k)$-vertex has at most
one $2$-neighbor. Hence the assumptions of Lemma~\ref{lem:11-neighbors} are satisfied and $G$ contains an edge of weight at most 13.
But this edge is nice and we get the contradiction with Invariant~\ref{inv} and the fact that $Q_e$ is empty.
\end{proof}

\begin{proposition}
\label{prop:nlogn}
The above algorithm can be implemented in $O(n\log n)$ time.
\end{proposition}

\begin{proof}
First we show that each recursive call takes only $O(\log n)$ amortized time. 
Checking adjacency in Step~\ref{step:contract} can be easily done in $O(\log n)$ time e.g.\ by storing the neighbors of each vertex in a balanced tree.
Then adding and removing edges can be done in $O(\log n)$ time. It remains to consider recoloring the graph after 
going back from the recursion. Recall from the proof of Lemma~\ref{lem:atmosttwo2vertices} that during the recoloring the algorithm
checks the colors of a bounded number of edges and also recolors a bounded number of edges.
Finding a free color can be done in constant time after a linear-time preprocessing (see Section 2.1.3 in~\cite{edge_col} for details).

The last unclear issue is verifying whether there is a path of given color, say $a$, between two vertices, say $x,y$.
W.l.o.g.\ we can assume that both $x$ and $y$ are incident with an edge colored $a$ (otherwise, immediately, the answer in negative), so
in fact, given two edges of the same color we want to check whether they are on the same path in the linear forest of color $a$.
Note that during recoloring an edge of color $a$ to $b$ (say), some path in the linear forest of $a$ is split into two paths (possibly one of length 0),
and some two paths (possibly empty) of the linear forest of color $b$ are connected to one path.
In other words, we need a data structure which maintains a linear forest that can be updated after adding or removing an edge and
processes connectivity queries of the form ``are the edges $e_1$ and $e_2$ on the same path?''. There are several solutions
to this problem generalized to forests with $O(\log n)$ time complexity (amortized) both for updates and queries -- e.g.\ link-cut trees of Sleator and Tarjan~\cite{TS} or ET-trees of Henzinger and King~\cite{ET-trees}. We note that in the case of {\em linear} forests this time complexity can be also achieved by using simply a balanced BST tree with efficient merge and split operations, like e.g.\ splay trees. All these data structures take
only linear space with respect to the size of the linear forest. 

Since in one recursive call we perform a bounded number of path queries, this takes only $O(\log n)$ amortized time. 
\end{proof}

\subsection{An algorithm for $\Delta \le 10$}
\label{sec:alg-bounded}

Now we sketch an algorithm which finds a partition of any planar graph of maximum degree $\Delta=O(1)$ into $k=\max\{\ceil{\tfrac{\Delta}{2}}, 5\}$ linear forests. Our algorithm uses all the reducible configurations described in Section~\ref{sec:structure}. Recall that they are of bounded size. 
Hence it is easy to check in constant time, whether a given vertex $v$ belongs to given configuration, since if this is the case, this configuration
is a subgraph of the graph induced of all vertices at some bounded distance from $v$ and because $\Delta=O(1)$ this subgraph has bounded size.
Our algorithm uses a queue of reducible configurations, initialized in linear time. Then configurations are added to the queue after modifying the graph.
Since each modification decreases the size of the graph, and causes appearance of a bounded number of configurations, the total number of configurations
is linear. After finding a configuration (by just removing it from the queue in constant time), shrinking the graph (usually by removing an edge), and going back from the recursive call, extending the coloring of the shrinked graph to the original graph takes $O(\log n)$ time, as described in the proof of Proposition~\ref{prop:nlogn}.

\begin{corollary}
Let $k=\max\{\ceil{\tfrac{\Delta}{2}}, 5\}$.
The above algorithm finds a $k$-linear coloring of any planar graph of maximum degree $\Delta$ in $O(n \log n)$ time.
\end{corollary}

%
%
%
%

\bibliographystyle{abbrv}




\end{document}